\documentclass{article}
\usepackage{amssymb}
\usepackage{graphicx}
\usepackage{amsmath}
\usepackage{amssymb}
\usepackage{amsfonts}
\usepackage{amscd}
\usepackage{amstext}
\usepackage{verbatim}
\usepackage{color}
\usepackage{multicol}
\usepackage{amsmath,amsthm,verbatim,amssymb,amsfonts,amscd, graphicx}
\usepackage{graphics}
\usepackage{tikz}
\usetikzlibrary{arrows,decorations.pathmorphing,backgrounds,positioning,fit,petri}
\usetikzlibrary{graphs,matrix}
\title{Functors induced by  Cauchy extension of C$^\ast$-algebras \\[0.3cm]}
\author{{ Kourosh Nourouzi$^1$\thanks{Corresponding author, e-mail: nourouzi@kntu.ac.ir} , \,\,\,Ali Reza$^{2}$
}\\[0.4cm]
{\em $^{1,2}$ Faculty of Mathematics, K. N. Toosi University of}\\
{\em Technology,}\\
{\em P.O. Box 16315-1618, Tehran, Iran.}\\}

\newtheorem{definition}{Definition}{\rm}
\newtheorem{corollary}{Corollary}

\newtheorem{theorem}{Theorem}
\newtheorem{Proposition}{Proposition}
\newtheorem{remark}{Remark}

\DeclareMathOperator{\ima}{Im}

\begin{document}
\maketitle

\begin{abstract}
In this paper we give three functors $\mathfrak{P}$, $[\cdot]_K$ and $\mathfrak{F}$ on the category of C$^\ast$-
algebras. The functor $\mathfrak{P}$ assigns to each C$^\ast$-algebra $\mathcal{A}$ a pre-C$^\ast$-algebra $\mathfrak{P}(\mathcal{A})$ with completion $[\mathcal{A}]_K$. The functor $[\cdot]_K$ assigns to each C$^\ast$-algebra $\mathcal{A}$ the Cauchy extension $[\mathcal{A}]_K$ of $\mathcal{A}$ by a non-unital C$^\ast$-algebra $\mathfrak{F}(\mathcal{A})$. Some properties of these functors are also given.  In particular, we show that the functors $[\cdot]_K$ and $\mathfrak{F}$ are exact and the functor $\mathfrak{P}$ is normal exact.
\end{abstract}
\def\thefootnote{ \ }
\footnotetext{{\em} $2010$ Mathematics Subject Classification. 46L05, 46M15}
\section{Introduction}
Given a complex C$^\ast$-algebra $\mathcal{A}$, the algebra $\mathcal{A}[[Z]]$  consists of all sequences $(a_n)_{n=0}^{\infty}$ in $\mathcal{A}$ with pointwise linear operations and Cauchy product
$$ ((a_n)_{n=0}^{\infty})((b_n)_{n=0}^{\infty})=(c_n)_{n=0}^{\infty},$$
where each $c_n=\sum_{k=0}^{n}a_kb_{n-k}$. It is natural to think of elements of $\mathcal{A}[[Z]]$ as the formal power series  in one variable of the form $\sum_{n=0}^{\infty}a_nZ^n$ with product
$$(\sum_{n=0}^{\infty}a_nZ^n)(\sum_{n=0}^{\infty}b_nZ^n)=\sum_{n=0}^{\infty}c_nZ^n,$$
 where $c_n$'s are as above. One may consider the complex subalgebra
$$\mathcal{A}[Z]=\{\sum_{n=0}^{\infty}a_nZ^n: \sum_{n=0}^{\infty}\|a_n\|<\infty\},$$
of $\mathcal{A}[[Z]]$. It is of interest to find a C$^\ast$-algebra via $\mathcal{A}[Z]$ to be an extension of $\mathcal{A}$.
Recall that an extension  $\mathcal{B}$ of $\mathcal{C}$ by $\mathcal{A}$ is a short exact sequence
\begin{equation}\label{equation01}
0\longrightarrow\mathcal{A}\overset{f}{\longrightarrow}\mathcal{B}\overset{g}{\longrightarrow}\mathcal{C}\longrightarrow 0.
\end{equation}
of C$^\ast$-algebras (see, e.g., \cite{5,14,1,15}). For any subset  $K$of $[-1,1]$ such that  $0$ is a limit point of $K$, we will define a pre-C$^\ast$-norm on $\mathcal{A}[Z]$. The completion of $\mathcal{A}[Z]$, denoted by $[\mathcal{A}]_K$, is an extension  of $\mathcal{A}$ (Proposition $\ref{proposition07}$ (iii)) which will be called the Cauchy extension of $\mathcal{A}$.

The outline of this work is as follows. In Section $\ref{section03}$ we introduce pre-C$^\ast$-algebra $\mathcal{A}[Z]$. In Proposition \ref{proposition05}, it is shown that $\mathcal{A}[Z]$ is not a C$^\ast$-algebra. Proposition \ref{proposition07} shows that the completion $[\mathcal{A}]_K$ of pre-C$^\ast$-algebra $\mathcal{A}[Z]$ is an extension of $\mathcal{A}$. We also introduce the functors $\mathfrak{P}$, $[\cdot]_K$ and $\mathfrak{F}$ on the category of C$^\ast$-algebras. The functor $\mathfrak{P}$ assigns to each C$^\ast$-algebra $\mathcal{A}$ a pre-C$^\ast$-algebra $\mathfrak{P}(\mathcal{A})=\mathcal{A}[Z]$. The functor $[\cdot]_K$ assigns to each C$^\ast$-algebra $\mathcal{A}$ an extension $[\mathcal{A}]_K$ of $\mathcal{A}$ by a non-unital C$^\ast$-algebra $\mathfrak{F}(\mathcal{A})$, where the C$^\ast$-algebra $\mathfrak{F}(\mathcal{A})$ is the completion of the ideal
$$\mathcal{A}_1=\{\sum_{n=0}^{\infty}a_nZ^n\in\mathcal{A}[Z]:a_0=0\}$$
 of $\mathcal{A}[Z]$. Some properties of functors $\mathfrak{P}$, $[\cdot]_K$ and $\mathfrak{F}$ are listed in Proposition $\ref{proposition08}$. In Section $\ref{section04}$ we show that the functors $[\cdot]_K$ and $\mathfrak{F}$ are exact. In Section $\ref{section05}$, using the notion of normal exact sequence of the normed spaces introduced by Yang \cite{2}, we prove that the functor $\mathfrak{P}$ is normal exact. More precisely, for any short exact sequence of C$^\ast$-algebra $(\ref{equation01})$ the corresponding short exact sequence
 $$0\longrightarrow\mathcal{A}[Z]\overset{\tilde{f}}{\longrightarrow}\mathcal{B}[Z]\overset{\tilde{g}}{\longrightarrow}\mathcal{C}[Z]\longrightarrow 0$$
 is a normal  exact sequence of pre-C$^\ast$-algebras. That is,  $\mathcal{B}(Z)/\ker  \tilde{g}\longrightarrow\mathcal{C}[Z]$ is an isometry. Among other results we also show that for any closed ideal $\mathcal{I}$ of a C$^\ast$-algebra $\mathcal{A}$, the pre-C$^\ast$-algebra $\mathcal{I}[Z]$ is a closed ideal of $\mathcal{A}[Z]$ (Proposition \ref{proposition08} (iii)) and the quotient $\mathcal{A}[Z]/\mathcal{I}[Z]$ is a pre-C$^\ast$-algebra (Theorem \ref{theorem03}) which is isometric $\ast$-isomorphic to $(\mathcal{A}/\mathcal{I})[Z]$ (Theorem \ref{theorem04}). Finally in Section \ref{section06}, we show that the Cauchy extension $[\mathcal{A}]_K$ of a C$^\ast$-algebra $\mathcal{A}$ can be considered as a C$^\ast$-subalgebra of $C_b(K,\mathcal{A})$, the C$^\ast$-algebra of all bounded continuous functions from $K$ to $\mathcal{A}$ (Theorem \ref{theorem05} (i)). In particular, if $K$ is compact, then  $[\mathcal{A}]_K$ is $\ast$-isomorphic to  $ C(K,\mathcal{A})$. We also give some other results in Theorem $\ref{theorem05}$. A minimax type result is given in Corollary \ref{corollary05}.

\section{Cauchy extension of C$^\ast$- algebras}\label{section03}

Let $\mathcal{A}$ be a complex Banach algebra and $\mathcal{A}[[Z]]$ be the complex algebra consisting of all formal power series in $\mathcal{A}$. If $\mathcal{A}$ has a unit, then an element $F=F(Z)=\sum_{n=0}^{\infty}a_nZ^n\in\mathcal{A}[[Z]]$ is invertible if and only if $a_0$ is an invertible element in $\mathcal{A}$. In particular, $1+Z^2$ is invertible in $\mathcal{A}[[Z]]$ and we have
\begin{equation}\label{equation02}
  (1+Z^2)(\sum_{n=0}^{\infty}(-1)^nZ^{2n})=(\sum_{n=0}^{\infty}(-1)^nZ^{2n})(1+Z^2)=1.
\end{equation}
The subalgebra
 $$\mathcal{A}[Z]=\{\sum_{n=0}^{\infty}a_nZ^n\in \mathcal{A}[[Z]]: \sum_{n=0}^{\infty}\|a_n\|<\infty\}$$
can be equipped with a norm as
  \begin{equation}\label{equation03}
    \|F\|=\sum_{n=0}^{\infty}\|a_n\|,
  \end{equation}
for all $F(Z)=\sum_{n=0}^{\infty}a_nZ^n\in \mathcal{A}[Z].$

\begin{Proposition}\label{proposition01}
Let $\mathcal{A}$ be a Banach algebra. Then $\mathcal{A}[Z]$ with the norm given in $(\ref{equation03})$ is a Banach algebra.

\noindent \begin{proof}
 To show that $\mathcal{A}[Z]$ is a Banach algebra, let $(F_k)=(\sum_{n=0}^{\infty}a_{kn}Z^n)$ be a sequence in $\mathcal{A}[Z]$ such that $\sum_{k=0}^{\infty}\|F_k\|<\infty$. Then

   $$\sum_{k=0}^{\infty}\sum_{n=0}^{\infty}\|a_{kn}\|=\sum_{n=0}^{\infty}\sum_{k=0}^{\infty}\|a_{kn}\|<\infty .$$
  Let $c_n=\sum_{k=0}^{\infty}a_{kn}$ and $F=\sum_{n=0}^{\infty}c_nZ^n$. Then $F\in A[Z]$. Let  $\varepsilon >0$ be given. There exists a positive integer $N$ such that    $\sum_{k=N+1}^{\infty}\sum_{n=0}^{\infty}\|a_{kn}\|<\varepsilon$. We have

\begin{center}
\begin{tabular}{lll}
$\|\sum_{k=0}^{N}F_k-F\|$&$=$&$\|\sum_{n=0}^{\infty}(\sum_{k=N+1}^{\infty}a_{kn})Z^n\|$\\[0.2cm]
&$=$&$\sum_{n=0}^{\infty}\|\sum_{k=N+1}^{\infty}a_{kn}\|$\\[0.2cm]
&$\leq$&$\sum_{n=0}^{\infty}\sum_{k=N+1}^{\infty}\|a_{kn}\| $\\[0.2cm]
&$=$&$\sum_{k=N+1}^{\infty}\sum_{n=0}^{\infty}\|a_{kn}\|$\\[0.2cm]
&$<$&$\varepsilon.$\\[0.2cm]
\end{tabular}
\end{center}
This completes the proof.
\end{proof}
\end{Proposition}
\begin{Proposition}\label{proposition02}
  Let $\mathcal{A}$ be a Banach algebra. If  $F(Z)=\sum_{n=0}^{\infty}a_nZ^n\in \mathcal{A}[Z]$, then $\sum_{n=0}^{N}a_nZ^n\rightarrow F(Z)$ as $N\rightarrow \infty$,

\noindent \begin{proof}
 Since
      $$\|F(Z)-\sum_{n=0}^{N}a_nZ^n\|=\|\sum_{n=N+1}^{\infty}a_nZ^n\|=\sum_{n=N+1}^{\infty}\|a_n\|,$$
    we get the desired limit.
  \end{proof}
\end{Proposition}

Now one can consider any element  $F(Z)=\sum_{n=0}^{\infty}a_nZ^n\in \mathcal{A}[Z]$ as a convergent series in $\mathcal{A}[Z]$.

If  $\mathcal{A}$ is a  C$^\ast$-algebra, we can define an involution $\ast$ in $\mathcal{A}[Z]$ by $F^\ast(Z)=\sum_{n=0}^{\infty}a_n^\ast Z^n$ for any $F(Z)\in\mathcal{A}[Z]$. In this case, $\mathcal{A}[Z]$ equipped with  this involution  and the norm given in $(\ref{equation03})$ is a $\ast$-Banach algebra.

\begin{Proposition}\label{proposition03}
  Let $\mathcal{A}$ be a  C$^\ast$-algebra.   There is no norm on involutive algebra $(\mathcal{A}[Z],\ast)$ which makes it a C$^\ast$-algebra. In particular, $(\mathcal{A}[Z], \ast)$ equipped with the norm given in $(\ref{equation03})$ is not a C$^\ast$-algebra.

  \noindent \begin{proof}
   We suppose on the contrary that there exists a norm $\|\cdot\|$ such that $(\mathcal{A}[Z],\ast,\|\cdot\|)$ is a $C^\ast$-algebra. Suppose that $\mathcal{A}$ is unital. By $(\ref{equation02})$ the element $1+Z^2$ is not invertible in $\mathcal{A}[Z]$.  This implies that $-1\in \sigma (Z^2)$ which is a contradiction. Now let $\mathcal{A}$ be non-unital and $a\in\mathcal{A}$ be self-adjoint with $\|a\|>1$. Applying  $(\ref{equation02})$ for $aZ$ we get that  $1+a^2Z^2$ is not invertible in $(\mathcal{A}\oplus \mathbb{C})(Z)$. That is, $-1\in \sigma (a^2Z^2)$ which is again  a contradiction.
  \end{proof}
\end{Proposition}

For a C$^\ast$-algebra $(\mathcal{A}, \|\cdot\|)$ if $F(Z)=\sum_{n=0}^{\infty}a_nZ^n\in \mathcal{A}[Z]$ and $-1\leq t\leq 1$ then
$$\sum_{n=0}^{\infty}\|a_nt^n\|\leq \sum_{n=0}^{\infty}\|a_n\|<\infty.$$
Hence $F(t)=\sum_{n=0}^{\infty}a_nt^n$ is norm-convergent in $\mathcal{A}$. %For example, the interval $[-1,1]$ is the largest one for which the series $\sum_{n=0}^{\infty}1/(n+1)^2Z^n\in\mathcal{A}[Z]$ is convergent on it.

For any $F(Z), G(Z)\in\mathcal{A}[Z]$ and $\lambda\in\mathbb{C}, t\in[-1,1]$ we have
\begin{equation}\label{equation04}
  (\lambda F(Z))(t)=\lambda F(t),
\end{equation}
\begin{equation}\label{equation05}
  (F(Z)+G(Z))(t)=F(t)+G(t),
\end{equation}
 \begin{equation}\label{equation06}
 (F(Z)G(Z))(t)=F(t)G(t).
\end{equation}
Note that the equalities $(\ref{equation04})$ and $(\ref{equation05})$ are clear and the proof of $(\ref{equation06})$ is similar to that of complex case (see \cite[p. 74]{4}).

The following proposition has a straightforward proof which is omitted here.
\begin{Proposition}\label{proposition04}
  Suppose that $K$ is a subset of $ [-1,1]$ such that $0$ is a limit point of $K$ and $(a_n)_{n=0}^\infty$ is a sequence 
  in  C$^\ast$-algebra $\mathcal{A}$. If 
  \begin{description}
  \item[(i)] $F(Z)=\sum_{n=0}^{\infty}a_nZ^n$ such that $\sum_{n=0}^{\infty}\|a_n\|<\infty$;
  \item[(ii)] $F(t)=0$ for any $t\in K$,
\end{description}
then $a_n=0$ for all $n$.
  \end{Proposition}

Hereafter, throughout the paper $K$ will denote a subset of $[-1,1]$ such that $0$ is a limit point of it.
%\begin{definition}\label{definition01}\rm{
%  Define $\|\cdot\|_K$ on $\mathcal{A}[Z]$ by
%  $$ \|F\|_K= \sup_{t\in K}\|\sum_{n=0}^{\infty}a_nt^n\|,$$
%    for all $F=F(Z)=\sum_{n=0}^{\infty}a_nZ^n\in \mathcal{A}[Z].$}
%    \end{definition}
\begin{Proposition}\label{proposition05}
The following statements hold:
\begin{description}
\item[(i)] The functional  $\|\cdot\|_K$  defined by 
  $$ \|F\|_K= \sup_{t\in K}\|\sum_{n=0}^{\infty}a_nt^n\|,$$
    for all $F=F(Z)=\sum_{n=0}^{\infty}a_nZ^n\in \mathcal{A}[Z],$ is a norm; 
  \item[(ii)] $(\mathcal{A}[Z], \ast, \|\cdot\|_K)$ is a pre-C$^\ast$-algebra but not a C$^\ast$-algebra;
  \item[(iii)] $\|F\|_K\leq \|F\|$ for all $F\in \mathcal{A}[Z]$;
  \item[(iv)] If $F(Z)=\sum_{n=0}^{\infty}a_nZ^n$, then $\sum_{n=0}^{N}a_nZ^n\rightarrow F(Z)$ as $N\to\infty$ in $\|\cdot\|_K$.
\end{description}

  \noindent \begin{proof}
    (i) From $(\ref{equation04})$, $(\ref{equation05})$, $(\ref{equation06})$ and Proposition  $\ref{proposition04}$ it is easily seen that $\|\cdot\|_K$ is a norm. (ii) By the definition of $\|\cdot\|_K$ we have the identity $\|F^\ast F\|_K=\|F\|_K^2$. Therefore $(\mathcal{A}[Z], \ast, \|\cdot\|_K)$ is a pre-C$^\ast$-algebra which by Proposition $\ref{proposition03}$ is not a C$^\ast$-algebra. (iii)  By the definition of $\|\cdot\|_K$ is clear.
  (iv)  The proof  follows from Proposition $\ref{proposition02}$ and Part (ii).
  \end{proof}
\end{Proposition}

We will call the completion $[\mathcal{A}]_K$ of pre-C$^\ast$-algebra $(\mathcal{A}, \ast, \|\cdot\|_K)$ the \textit{K}-Cauchy or simply the Cauchy extension of $\mathcal{A}$. It is clear that $[\mathcal{A}]_K$ is a C$^\ast$-algebra. %Also if $K'\subseteq  K$ such that $0$ is a limit point of $K'$, then it is clear that $\|%\cdot\|_{K'}\leq\|\cdot\|_K$. In particular, if $J=[-1,1]$ then for any $K$ we have $\|\cdot\|_K\leq\|\cdot\|%_J$.
\begin{Proposition}\label{proposition06}
Let $\mathcal{A}$ be a C$^\ast$-algebra. The following  hold:
  \begin{description}
    \item[(i)]
 If $\mathfrak{I}$ is an ideal of $\mathcal{A}[Z]$, then the completion  $\hat{\mathfrak{I}}$ of $(\mathfrak{I}, \|\cdot\|_K)$ is a closed ideal of $[\mathcal{A}]_K$;
    \item[(ii)] If $\mathcal{I}$ is a closed ideal of $\mathcal{A}$, then $[\mathcal{I}]_K$ is a closed ideal of $[\mathcal{A}]_K$.

  \end{description}

  \noindent \begin{proof}
 (i) Let $\mathfrak{I}$ be an ideal of $\mathcal{A}[Z]$. Then the completion  $\hat{\mathfrak{I}}$ of $(\mathfrak{I}, \|\cdot\|_K)$ is a closed ideal of $[\mathcal{A}]_K$. Choose any element $F\in\hat{\mathfrak{I}}$ and $G\in [\mathcal{A}]_K$. Let $(F_n)$ and $(G_k)$ be two sequences in $\mathfrak{I}$ and $\mathcal{A}[Z]$ respectively converging to $F\in\hat{\mathfrak{I}}$ and $G\in[\mathcal{A}]_K$. For any $k,n\geq 1$ we have $F_nG_k, G_kF_n\in\mathfrak{I}$. This implies that $FG_K, G_kF\in\hat{\mathfrak{I}}$, for all $k\geq 1$ and so $FG, GF\in\hat{\mathfrak{I}}$. That is $\hat{\mathfrak{I}}$ is a closed ideal of $[\mathcal{A}]_K$.\\
    (ii) Consider $F\in \mathcal{I}[Z]$ and $G\in \mathcal{A}[Z]$. It is clear that $FG,GF\in \mathcal{I}[Z]$, i.e., $\mathcal{I}[Z]$ is an ideal of  $\mathcal{A}[Z]$. Now Part (i) implies that $\widehat{(\mathcal{I}[Z])}=[\mathcal{I}]_K$ is a closed ideal of $[\mathcal{A}]_K$.

\end{proof}
\end{Proposition}

For a C$^\ast$-algebra $\mathcal{A}$ define
$$\mathcal{A}_0=\{F(Z)=\sum_{n=0}^{\infty}a_nZ^n\in \mathcal{A}[Z]:a_n=0 \ \ \ for \ \ \  n>0\},$$

  $$\mathcal{A}_1=\{F(Z)=\sum_{n=0}^{\infty}a_nZ^n \in \mathcal{A}[Z]:a_0=0\}.$$
Denote the completion of $\mathcal{A}_1$ by $\hat{\mathcal{A}}_1$. It is clear that $\mathcal{A}_1$ is an ideal of $\mathcal{A}[Z]$ and by Proposition $\ref{proposition06}$,
$\hat{\mathcal{A}}_1$ is a closed ideal of $[\mathcal{A}]_K$. Hence if $\mathcal{A}\neq 0$, then $[\mathcal{A}]_K$ has a proper closed ideal $\hat{\mathcal{A}}_1$. Consequently no simple C$^\ast$-algebra is a Cauchy extension of some
 C$^\ast$-algebra.
 It is worth mentioning that there is no ideal $\mathcal{I}$ of $\mathcal{A}$ such that $\mathcal{I}[Z]=\mathcal{A}_1$.  Since $\mathcal{A}_0$ is naturally $\ast$-isomorphic to $\mathcal{A}$ we always use $\mathcal{A}$ instead of $\mathcal{A}_0$ as a subalgebra of $\mathcal{A}[Z]$.

Suppose that $\mathcal{A}, \mathcal{B}, \mathcal{E}$ are C$^\ast$-algebras such that $\mathcal{B}$ is an ideal of $\mathcal{E}$. It is said to be $\mathcal{E}$  an extension of $\mathcal{A}$ by $\mathcal{B}$ if there is a short exact sequence
$$0\longrightarrow \mathcal{B}\overset{i}{\longrightarrow} \mathcal{E}\overset{p}{\longrightarrow} \mathcal{A}\longrightarrow 0,$$
where $i(\mathcal{B})=\ker p$ and $i, p$ are injective and surjective $\ast$-homomorphisms respectively (see, e.g.,\cite{5}).

\begin{Proposition}\label{proposition07}
Let $\mathcal{A}$ be a C$^\ast$-algebra. The following statements hold:
  \begin{description}
    \item[(i)] Every element $F$ of $[\mathcal{A}]_K$ has a unique representation $F=a+G$, where $a\in \mathcal{A}$ and $G\in \hat{\mathcal{A}_1}$;
    \item[(ii)] $\|a\|_K=\|a\|\leq \|a+G\|_K$, for all $a\in\mathcal{A}$ and $G\in\hat{\mathcal{A}}_1$;
    \item[(iii)] $[\mathcal{A}]_K$ is an extension of $\mathcal{A}$ by $\hat{\mathcal{A}}_1$;
     \item[(iv)] $\hat{\mathcal{A}}_1$ is not unital as a $C^\ast$-subalgebra of $[\mathcal{A}]_K$.
  \end{description}

  \noindent \begin{proof} (i) Let $(F_k)$ be a Cauchy sequence in $(\mathcal{A}[Z], \|\cdot\|_K)$, where 
  $$F_k=\sum_{n=0}^{\infty}a_{kn}Z^n\in\mathcal{A}[Z].$$
   Let $\varepsilon >0$ be given. Then $\|F_k-F_{k'}\|_K<\varepsilon$ for sufficiently large $k, k'$ . Suppose that $(t_m)$ is a sequence in $K$ such that $t_m\longrightarrow 0$ as $m\longrightarrow\infty$. By the definition of $\|\cdot\|_K$ we have

    $$\|a_{k0}-a_{k'0}\|=\lim_{m\to \infty}\|\sum_{n=0}^{\infty}(a_{kn}-a_{k'n})t^n_m\|\leq \sup_{t\in K}\|\sum_{n=0}^{\infty}(a_{kn}-a_{k'n})t^n\|$$
    $$=\| F_k-F_{k'}\|_K<\varepsilon ,$$
    for sufficiently large $k, k'$. Furthermore
    $$\sup _{t\in K} \|\sum_{n=1}^{\infty}(a_{kn}-a_{k'n})t^n\|<2\varepsilon .$$
Therefore the sequences $(a_{k0})$ and $(\sum_{n=1}^{\infty}a_{kn}Z^n)$ are Cauchy  in $\mathcal{A}$ and $\mathcal{A}_1$, respectively. For $F\in [\mathcal{A}]_K$, let $F=\lim_{k\to\infty}F_k$, where $F_k=\sum^{\infty}_{n=0}a_{kn}Z^n\in\mathcal{A}[Z]$. Then $F=a+G$, where $a_{k0}\longrightarrow a\in\mathcal{A}$ and $\sum_{n=1}^{\infty}a_{kn}Z^n\longrightarrow G\in \hat{\mathcal{A}}_1$ as $k\to \infty$. Since $\hat{\mathcal{A}}_1\cap\mathcal{A}=0$, then this representation is unique. Hence $[\mathcal{A}]_K$ is the internal direct sum of subspaces $\mathcal{A}$ and $\hat{\mathcal{A}_1}$, i.e., $[\mathcal{A}]_K=\mathcal{\mathcal{A}}\oplus \hat{\mathcal{A}_1}$.\\
    (ii) Note that if $a\in\mathcal{A}$ and $G=\lim_{k\to\infty}\sum_{n=1}^{\infty}a_{kn}Z^{n}\in\hat{\mathcal{A}}_1$, then
    $$\|a+\sum_{n=1}^{\infty}a_{kn}Z^n\|_K=\sup_{t\in K}\|a+\sum_{n=1}^{\infty}a_{kn}t^n\|,$$
   for all $k\geq 1$. A similar method to that used in Part (i) implies that  $\|a\|\leq \|a+\sum_{n=1}^{\infty}a_{kn}Z^n\|$, for all $k\geq 1$. Therefore $\|a\|\leq \|a+G\|$, for all $a\in\mathcal{A}$ and $G\in \hat{\mathcal{A}_1}$.\\
    (iii)   Define $p_\mathcal{A}:[\mathcal{A}]_K\longrightarrow \mathcal{A}$ by $p_\mathcal{A}(a+G)=a$, for all $a\in\mathcal{A}$ and $G\in\hat{\mathcal{A}}_1$. It is easily seen that $p_\mathcal{A}$ is a surjective $\ast$-homomorphism and $\ker p_\mathcal{A}=\hat{\mathcal{A}}_1$. Therefore we have the short exact sequence
      \begin{equation}\label{equation08}
  0\rightarrow\hat{\mathcal{A}_1}\overset{i}{\hookrightarrow}[\mathcal{A}]_K\overset{p_A}{\rightarrow}\mathcal{A}\rightarrow 0.
  \end{equation}
 This shows that $[\mathcal{A}]_K$ is an extension of $\mathcal{A}$ by $\hat{\mathcal{A}}_1$.\\
 (iv) Suppose on the contrary that $\hat{\mathcal{A}_1}$ is unital with unit $U(Z)$. Since $aZU(Z)=aZ$ for all $a\in A$, we have $taU(t)=ta$ for any $t\in K$ and $a\in A$. This implies that $aU(t)=a$ for all $t\neq 0$ and therefore $\lim_{t\to 0}U(t)\neq 0$, which is a contradiction.

  \end{proof}
\end{Proposition}

\begin{remark}\label{definition02}\rm{
 Each  $\ast$-homomorphism $f:\mathcal{A}\longrightarrow \mathcal{B}$ of C$^\ast$-algebras induces a  $\ast$-homomorphism $\tilde{f}:\mathcal{A}[Z]\longrightarrow \mathcal{B}[Z]$  between pre-C$^\ast$-algebras $\mathcal{A}[Z]$ and $\mathcal{B}[Z]$ by
 \begin{equation}\label{equation09}
 \tilde{f}(\sum_{n=0}^{\infty}a_nZ^n)=\sum_{n=0}^{\infty}f(a_n)Z^n,
 \end{equation}
where $\sum_{n=0}^{\infty}a_nZ^n\in\mathcal{A}[Z]$.}
\end{remark}

\begin{remark}\label{remark02}\rm{If we define $\mathfrak{P}(\mathcal{A})=\mathcal{A}[Z]$ for any C$^\ast$-algebra $\mathcal{A}$ and $\mathfrak{P}(f)=\tilde{f}$, for any $\ast$-homomorphism  $f:\mathcal{A}\longrightarrow\mathcal{B}$ of C$^\ast$-algebras, then $\mathfrak{P}$ is a functor from the category of C$^\ast$-algebras to the category of pre-C$^\ast$-algebras.
Each $\ast$-homomorphism $\tilde{f}:\mathcal{A}[Z]\longrightarrow \mathcal{B}[Z]$ defined by $(\ref{equation09})$ induces a $\ast$-homomorphism $\hat{f}:[\mathcal{A}]_K\longrightarrow [\mathcal{B}]_K$. It is easy to see that $[\cdot]_K$ is a functor from the category of C$^\ast$-algebras to itself as $[f]_K=\hat{f}$. Now defining  $\mathfrak{F}(\mathcal{A})=\hat{\mathcal{A}}_1$ and $\mathfrak{F}(\mathcal{A}\overset{f}{\longrightarrow}\mathcal{B})=\hat{f}|_{\hat{\mathcal{A}}_1}:\hat{\mathcal{A}}_1\longrightarrow\hat{\mathcal{B}}_1$, for C$^\ast$-algebras $\mathcal{A}$, $\mathcal{B}$ and $\ast$-homomorphism $f$, we get a   functor on the category of C$^\ast$-algebras which  assigns, by Proposition $\ref{proposition07}$ (iv), to any C$^\ast$-algebra a non-unital C$^\ast$-algebra.}
\end{remark}

By $\mathcal{A}\cong \mathcal{B}$ we mean that the C$^\ast$-algebras $\mathcal{A}$ and $\mathcal{B}$ are $\ast$-isomorphic.

\begin{Proposition}\label{proposition08}
  Let $f:\mathcal{A}\longrightarrow \mathcal{B}$ be a $\ast$-homomorphism of C$^\ast$-algebras. Then

  \begin{description}
    \item[(i)] $\tilde{f}$ is a contraction;
    \item[(ii)] $\tilde{f}$ and $\hat{f}$ are isometries provided that $f$ is an isometry;
    \item[(iii)]  $\tilde{f}$ is surjective provided that $f$ is surjective;
    \item[(iv)] If $f$ is a $\ast$-isomorphism, then both $\tilde{f}$ and $\hat{f}$ are $\ast$-isomorphisms;
    \item[(v)] $\ker\tilde{f}=(\ker f)[Z]$;
    \item[(vi)]$\ima\tilde{f}=(\ima f)[Z]$;
    \item[(vii)] If $\mathcal{I}$ is a closed ideal of $\mathcal{A}$, then $\mathcal{I}[Z]$ is a closed ideal of $(\mathcal{A}[Z], \|\cdot\|_K)$. In particular,
    $$0\longrightarrow\mathcal{I}[Z]\hookrightarrow\mathcal{A}[Z]\overset{p'}{\longrightarrow}\mathcal{A}[Z]/\mathcal{I}[Z]\longrightarrow 0$$
and
$$0\longrightarrow\mathcal{I}[Z]\hookrightarrow\mathcal{A}[Z]\overset{\tilde{p}}{\longrightarrow}(\mathcal{A}/\mathcal{I})[Z]\longrightarrow 0$$
are short exact sequences;
    \item[(viii)]$[\mathcal{A}\oplus\mathcal{B}]_K\cong[\mathcal{A}]_K\oplus[\mathcal{B}]_K$;
  \item[(ix)]$\widehat{(\mathcal{A}\oplus\mathcal{B})}_1\cong\hat{\mathcal{A}}_1\oplus\hat{\mathcal{B}}_1$.
  \end{description}
  \noindent \begin{proof}The proof of (iv) follows from (ii) and (iii). The proofs of (v) and (vi) are straightforward and the proof of (ix) is similar to Part (viii). We prove the others.
   \begin{description}
    \item[(i)] suppose that $F(Z)=\sum_{n=0}^{\infty}a_nZ^n\in \mathcal{A}[Z]$. Then
\begin{center}
\begin{tabular}{lll}
$\|\tilde{f}(F)\|_K$&$=$&$\|\sum_{n=0}^{\infty}f(a_n)Z^n\|_K$\\[0.2cm]
&$=$&$\sup_{t\in K}\|\sum_{n=0}^{\infty}f(a_n)t^n\|$\\[0.2cm]
&$=$&$\sup_{t\in K}\|f(\sum_{n=0}^{\infty}a_nt^n)\|$\\[0.2cm]
&$\leq$&$\sup_{t\in K}\|\sum_{n=0}^{\infty}a_nt^n\|$\\[0.2cm]
&$=$&$\|F\|_K.$\\[0.2cm]
\end{tabular}
\end{center}
    \item[(ii)] If $f$ is an isometry, then the proof of (i) shows that $\|\tilde{f}(F)\|_K=\|F\|_K$, for all $F\in\mathcal{A}[Z]$.  That is $\tilde{f}$ and consequently $\hat{f}$ is an isometry.
    \item[(iii)] Let $f$ be surjective and $G=\sum_{n=0}^{\infty}b_nZ^n\in \mathcal{B}[Z]$. For any integer $n\geq 0$, there exists $a_n\in \mathcal{A}$ such that $b_n=f(a_n)$. For any integer $n\geq 0$ there exists $a_n'\in\ker f$ such that
      \begin{equation}\label{equation10}
    \|a_n+a'_n\|\leq \|a_n+\ker f\|+2^{-n}.
  \end{equation}
  Since $\mathcal{A}/\ker f\cong\mathcal{B}$, we have
       \begin{equation}\label{equation11}
\|a_n+\ker f\|=\|f(a_n)\|=\|b_n\|.
  \end{equation}
    Define $a_n''=a_n+a_n'$, for all $n\geq 0$. Now we see from $(\ref{equation10})$ and $(\ref{equation11})$ that $F(Z)=\sum_{n=0}^{\infty}a_n''Z^n\in\mathcal{A}[Z]$ and $f(a_n'')=b_n$, for each $n\geq 0$, and therefore $\tilde{f}(F)=G$.
    \item[(vii)]Exactness of first diagram is clear. Part (iii) shows that $\mathcal{A}[Z]\overset{\tilde{p}}{\rightarrow}(\mathcal{A}/\mathcal{I})[Z]$ induces by the projection $\mathcal{A}\overset{p}{\longrightarrow}\mathcal{A}/\mathcal{I}$ is surjective. By (v) $\ker\tilde{p}=\mathcal{I}[Z]$ is a closed ideal of $\mathcal{A}[Z]$. This completes the proof.
 \item[(viii)] It is easily seen that
 $$T:\mathcal{A}[Z]\oplus\mathcal{B}[Z]\longrightarrow(\mathcal{A}\oplus\mathcal{B})[Z]$$
 defined by
 $$T(\sum_{n=0}^{\infty}a_nZ^n, \sum_{n=0}^{\infty}b_nZ^n)=\sum_{n=0}^{\infty}(a_n, b_n)Z^n,$$
 for all $\sum_{n=0}^{\infty}a_nZ^n\in\mathcal{A}[Z]$ and $\sum_{n=0}^{\infty}b_nZ^n\in\mathcal{B}[Z]$, is a $\ast$-isomorphism.

    \end{description}
  \end{proof}
\end{Proposition}

\section{Exactness of the functor $[\cdot]_K$ }\label{section04}
In this section we show that $[\cdot]_K$ is an exact functor. We first recall some definitions of the 
category theory \cite{7}.

 Recall that a map $X\overset{f}{\longrightarrow}Y$ in a category $\mathfrak{C}$ is called an epimorphism if for all maps $Y\overset{g}{\longrightarrow}Z$ and $Y\overset{h}{\longrightarrow}Z$ in $\mathfrak{C}$ with $g\circ f=h\circ f$, we have $g=h$. In the category of C$^\ast$-algebras,  a $\ast$-homomorphism $f:\mathcal{A}\rightarrow\mathcal{B}$ is an epimorphism if and only if it is surjective \cite{6}.% We use this fact %to the Theorem $\ref{theorem01}$.

Suppose that $X\overset{f}{\longrightarrow}Y$ is a map in a category $\mathfrak{C}$ with zero object. A map $Z\overset{j}{\longrightarrow}X$ is a kernel of $f$ if $f\circ j=0$ and for any map $Z'\overset{g}{\longrightarrow}X$ in $\mathfrak{C}$ such that $f\circ g=0$, there exists a unique map $Z'\overset{h}{\longrightarrow}Z$ such that $j\circ h=g$. For example, if $\mathcal{A}\overset{f}{\longrightarrow}\mathcal{B}$ is a $\ast$-homomorphism of C$^\ast$-algebras, then the inclusion $\ker f\hookrightarrow\mathcal{A}$ is a kernel of $f$.

\begin{theorem}\label{theorem01}
 The functor $[\cdot]_K$ is exact.

\noindent \begin{proof}
Suppose that
$$0\longrightarrow\mathcal{A}\overset{f}{\longrightarrow}\mathcal{B}\overset{g}{\longrightarrow}\mathcal{C}\longrightarrow 0$$
is a short exact sequence of C$^\ast$-algebras. We must show that
\begin{equation}\label{equation12}
0\longrightarrow[\mathcal{A}]_K\overset{\hat{f}}{\longrightarrow}[\mathcal{B}]_K\overset{\hat{g}}{\longrightarrow}[\mathcal{C}]_K\longrightarrow 0
\end{equation}
is a short exact sequence of C$^\ast$-algebras. We first  show that if $f:\mathcal{A}\rightarrow\mathcal{B}$ is a surjective $\ast$-homomorphism of C$^\ast$-algebras, then $\hat{f}:[\mathcal{A}]_K\rightarrow[\mathcal{B}]_K$ is also a surjective $\ast$-homomorphism of C$^\ast$-algebras. To do this suppose that $[\mathcal{B}]_K\overset{h}{\rightarrow}\mathcal{C}$ and $[\mathcal{B}]_K\overset{g}{\rightarrow}\mathcal{C}$ are $\ast$-homomorphism of C$^\ast$-algebras such that $h\circ \hat{f}=g\circ \hat{f}$. From Proposition $\ref{proposition08}$ (iii), we have $\hat{f}(\mathcal{A}[Z])=\mathcal{B}[Z]$. So for any $G(Z)\in \mathcal{B}[Z]$ there exists an element $F(Z)\in\mathcal{A}[Z]$ such that
  \begin{center}
\begin{tabular}{lll}
$h(G(Z))$&$=$&$h(\hat{f}(F(Z)))$\\[0.2cm]
&$=$&$(g\circ \hat{f})(F(Z))$\\[0.2cm]
&$=$&$g(G(Z))$.\\[0.2cm]
\end{tabular}
\end{center}
    This implies that $h|_{\mathcal{B}[Z]}=g|_{\mathcal{B}[Z]}$ and therefore $g=h$. Hence $\hat{f}$ is  an epimorphism and consequently is surjective by \cite{6}.

Now we show that if $\mathcal{A}\overset{f}{\longrightarrow}\mathcal{B}$ is a $\ast$-homomorphism of C$^\ast$-algebras, then $\ker\hat{f}=[\ker f]_K$.
To prove this, suppose that $\mathcal{C}\overset{g}{\longrightarrow}[\mathcal{A}]_K$ is a $\ast$-homomorphism of C$^\ast$-algebras such that $\hat{f}\circ g=0$. If $\hat{g}_\mathcal{C}=\hat{p}_\mathcal{A}\circ \hat{g}$ and $\mathcal{C}\overset{i_\mathcal{C}}{\hookrightarrow}[\mathcal{C}]_K$ is the injection, then $\hat{g}_\mathcal{C}\circ i_\mathcal{C}=g$. Since $\ker f\overset{j}{\hookrightarrow}\mathcal{A}$ is a kernel of $f$, there exists a unique $\ast$-homomorphism $\mathcal{C}\overset{h}{\longrightarrow}\ker f$ such that the  diagram

\begin{center}
\begin{tikzpicture}[
scale=0.5,
ca/.style={shorten <=28pt},
cb/.style={shorten <=1.5pt},]
\matrix (m) [matrix of math nodes,row sep=3em,column sep=3em,minimum width=2em]
 {
     \ker f  & \mathcal{A} & \mathcal{B} \\
     & \mathcal{C} \\};
   \path[-stealth]
    (m-2-2) edge node [left] {$h$} (m-1-1)
     (m-1-1.east|-m-1-2) edge [cb] (m-1-2)
     (m-1-2) edge node [above] {$f$} (m-1-3)
     (m-2-2) edge node [right] {$g_{\mathcal{C}}$} (m-1-2);
     \path[-right hook]
     (m-1-2) edge [ca] (m-1-1.east|-m-1-2)   ;
\end{tikzpicture}
\end{center}
is commutative. Since $[\cdot]_K$ is a functor we get the commutative diagram

\begin{center}
\begin{tikzpicture}[
scale=0.5,
ca/.style={shorten <=28pt},
cb/.style={shorten <=1.5pt},]
\matrix (m) [matrix of math nodes,row sep=3em,column sep=3em,minimum width=2em]
 {
     [\ker f]_K  &[] [\mathcal{A}]_K &[] [\mathcal{B}]_K \\
     & [][\mathcal{C}]_K \\};
   \path[-stealth]
    (m-2-2) edge node [left] {$\hat{h}$} (m-1-1)
     (m-1-1.east|-m-1-2) edge [cb] node [above] {$j$}(m-1-2)
     (m-1-2) edge node [above] {$\hat{f}$} (m-1-3)
     (m-2-2) edge node [right] {$\hat{g}_{\mathcal{C}}$} (m-1-2);
     \path[-right hook]
     (m-1-2) edge [ca] (m-1-1.east|-m-1-2)   ;
\end{tikzpicture}
\end{center}
Putting $h'=\hat{h}\circ i_\mathcal{C}$ we get $j\circ h'=g$, since $j\circ\hat{h}=\hat{g}_\mathcal{C}$. Now we show that $h'$ is unique. Suppose that there is a $\ast$-homomorphism $\mathcal{C}\overset{k}{\longrightarrow} [\ker f]_K$  such that $j\circ k=g=j\circ h'$. Since $j$ is an injection, then $k=h'$, which proves the uniqueness of $h'$. It is clear that $\ker\hat{f}=[\ker f]_K$. Now the Parts (ii), (v) and (vi) of Proposition $\ref{proposition08}$ imply that $(\ref{equation12})$ is a short exact sequence of C$^\ast$-algebras, or equivalently $[\cdot]_K$ is an exact functor. The  diagram

\begin{center}
\begin{tikzpicture}[
scale=0.5,
ca/.style={shorten <=38pt},
cb/.style={shorten <=2pt},
cc/.style={shorten >=5pt},
cd/.style={shorten <=98pt}
]

\matrix (m) [matrix of math nodes,row sep=4em,column sep=4em,minimum width=2em]
{
[\ker f]_K& &[][\mathcal{A}]_K& &[][\mathcal{B}]_K \\
&[][\mathcal{C}]_K&&& \\
&\mathcal{C}&&& \\ };
\path[-stealth]
(m-1-1) edge [cb] node [above] {$j$} (m-1-3)
(m-3-2) edge [cb] node [left] {$i_\mathcal{C}$}(m-2-2)
(m-3-2) edge node [below left] {$k$} (m-1-1)
(m-2-2) edge node [above left] {$\hat{g}_\mathcal{C}$} (m-1-3)
(m-3-2) edge node [below right] {$g$}(m-1-3)
(m-1-3) edge node [above] {$\hat{f}$} (m-1-5)
(m-2-2) edge [cc] node [above right] {$\hat{h}$}(m-1-1)
;
\path[-right hook]
(m-2-2) edge [ca] (m-3-2);
\path[-right hook]
(m-1-3) edge [cd] (m-1-1) ;
\end{tikzpicture}
\end{center}
shows the detials above.
\end{proof}
\end{theorem}
\begin{corollary}\label{corollary01}
If $\mathcal{I}$ is a closed ideal of a C$^\ast$-algebra $\mathcal{A}$, then $[\mathcal{A}/\mathcal{I}]_K\cong[\mathcal{A}]_K/[\mathcal{I}]_K.$

\noindent \begin{proof}
By Theorem $\ref{theorem01}$, the short exact sequence
$$0\longrightarrow \mathcal{I}\hookrightarrow\mathcal{A}\longrightarrow\mathcal{A}/\mathcal{I}\longrightarrow 0 $$
induces the short exact sequence
$$0\longrightarrow [\mathcal{I}]_K\hookrightarrow[\mathcal{A}]_K\longrightarrow[\mathcal{A}/\mathcal{I}]_K\longrightarrow 0 $$
which implies that $[\mathcal{A}]_K/[\mathcal{I}]_K\cong [\mathcal{A}/\mathcal{I}]_K$.
\end{proof}
\end{corollary}

In the following corollary we use $3\times 3$ lemma in homological algebra for the C$^\ast$-algebras as complex vector spaces (see, e.g., \cite{8}).
\begin{corollary}\label{corollary02}
If
$$0\longrightarrow \mathcal{A}\overset{f}{\longrightarrow}\mathcal{B}\overset{g}{\longrightarrow}\mathcal{C}\longrightarrow 0$$
is a short exact sequence of C$^\ast$-algebras, then
$$0\longrightarrow \hat{\mathcal{A}}_1\overset{\hat{f}|_{\hat{\mathcal{A}}_1}}{\longrightarrow}\hat{\mathcal{B}}_1\overset{\hat{g}|_{\hat{\mathcal{B}}_1}}{\longrightarrow}\hat{\mathcal{C}}_1\longrightarrow 0$$
is also a short exact sequence of C$^\ast$-algebras, i.e., $\mathfrak{F}$ is an exact functor (see Remark $\ref{remark02}$). Furthermore, if $\mathcal{I}$ is a closed ideal of $\mathcal{A}$, then $(\widehat{\mathcal{A}/\mathcal{I}})_1\cong\hat{\mathcal{A}}_1/\hat{\mathcal{I}_1}$

\noindent \begin{proof}
In the commutative diagram

\begin{center}
\begin{tikzpicture}[
scale=0.5,
ca/.style={shorten <=28.3pt},
cb/.style={shorten <=2pt},]
\matrix (m) [matrix of math nodes,row sep=4em,column sep=4em,minimum width=2em]
 {
     &0&0&0& \\
     0&\hat{\mathcal{A}}_1 &\hat{\mathcal{B}}_1&\hat{\mathcal{C}}_1&0 \\
     0& [][\mathcal{A}]_K &[][\mathcal{B}]_K&[][\mathcal{C}]_K&0 \\
     0&\mathcal{A} &\mathcal{B}&\mathcal{C}&0 \\
     &0&0&0& \\ };
   \path[-stealth]
   (m-1-2) edge (m-2-2)
   (m-1-3) edge (m-2-3)
   (m-1-4) edge (m-2-4)
   (m-3-1) edge (m-3-2)
   (m-3-2) edge node [above] {$\hat{f}$} (m-3-3)
   (m-3-3) edge node [above] {$\hat{g}$}(m-3-4)
   (m-3-2) edge node [right] {$p_\mathcal{A}$} (m-4-2)
   (m-3-3) edge node [right] {$p_\mathcal{B}$} (m-4-3)
   (m-3-4) edge node [right] {$p_\mathcal{C}$} (m-4-4)
   (m-3-4) edge (m-3-5)
   (m-4-1) edge (m-4-2)
   (m-4-2) edge node [above] {$f$} (m-4-3)
   (m-4-3) edge node [above] {$g$} (m-4-4)
   (m-4-4) edge (m-4-5)
   (m-4-2) edge (m-5-2)
   (m-4-3) edge (m-5-3)
   (m-4-4) edge (m-5-4)
   (m-2-1) edge [->,dotted] (m-2-2)
   (m-3-2) edge [ca,-right hook]  (m-2-2)
   (m-2-2) edge [cb] (m-3-2)
   (m-3-3) edge [ca,-right hook]  (m-2-3)
   (m-2-3) edge [cb] (m-3-3)
   (m-3-4) edge [ca,-right hook]  (m-2-4)
   (m-2-4) edge [cb] (m-3-4)
   (m-2-4) edge [->,dotted] (m-2-5)
   (m-2-2) edge node [above] {$\hat{f}|_{\hat{\mathcal{A}}_1}$} (m-2-3)
   (m-2-3) edge node [above] {$\hat{g}|_{\hat{\mathcal{B}}_1}$} (m-2-4)
   ;
\end{tikzpicture}
\end{center}
the middle row is exact by Theorem $\ref{theorem01}$ and all columns are exact by $(\ref{equation08})$. Now $3\times 3$ Lemma \cite{8} shows that the top row is also exact. By a similar argument as in Corollary, $\ref{corollary01}$ we get $$(\widehat{\mathcal{A}/\mathcal{I}})_1\cong\hat{\mathcal{A}}_1/ \hat{\mathcal{I}_1}.$$

\end{proof}
\end{corollary}

Recall that an ideal $\mathcal{I}$ of a C$^\ast$-algebra $\mathcal{A}$ is called modular if there is an element $u\in\mathcal{A}$ such that $ua-a, au-a\in\mathcal{A}$, for all element $a\in\mathcal{A}$. Note that $\mathcal{I}$ is modular if and only if ${\mathcal{A}}/{\mathcal{I}}$ is unital \cite{9}.

\begin{corollary}\label{corollary03}
Let $\mathcal{I}$ be a closed ideal of a C$^\ast$-algebra $\mathcal{A}$. Then $\mathcal{I}$ is a modular ideal of $\mathcal{A}$ if and only if $[\mathcal{I}]_K$ is a modular ideal of $[\mathcal{A}]_K$.

\noindent \begin{proof}
We first  show that a C$^\ast$-algebra $\mathcal{B}$ is unital if and only if $[\mathcal{B}]_K$ is unital. It can be easily seen that if $\mathcal{B}$ is unital, then $[\mathcal{B}]_K$ is also unital. Now, by Proposition $\ref{proposition07}$ (i) suppose that $[\mathcal{B}]_K$ is unital with unit $a+G$ for some $a\in\mathcal{B}$ and $G\in\hat{\mathcal{B}}_1$. Consider an arbitrary element $b+F\in [\mathcal{B}]_K$ with $b\in\mathcal{B}$ and $F\in\hat{\mathcal{B}}_1$. Then $(b+F)(a+G)=b+F$ or equivalently $ba+FG+Fa+bG=b+F$. It follows that $ba-b=H$, for some $H\in\hat{\mathcal{B}}_1$. Since $\mathcal{B}\cap\hat{\mathcal{B}}_1=0$, then  $ba=b$. Similarly $ab=b$. This shows that $a$ is the unit of $\mathcal{B}$. Now let $\mathcal{I}$ be a closed ideal of $\mathcal{A}$. Then by Corollary $\ref{corollary01}$, $\mathcal{I}$ is modular if and only if $[{\mathcal{A}}/{\mathcal{I}}]_K\cong {[\mathcal{A}]_K}/{[\mathcal{I}]_K}$ is unital. Hence $\mathcal{I}$ is modular if  and only if $[\mathcal{I}]_K$ is modular.
\end{proof}
\end{corollary}

\section{Normal exactness of the functor $\mathfrak{P}$}\label{section05}
Suppose that $\mathcal{A}$ is a C$^\ast$-algebra and $\mathcal{I}$ is a closed ideal of $\mathcal{A}$. It follows from Proposition $\ref{proposition08}$ (vii) that ${\mathcal{A}[Z]}/{\mathcal{I}[Z]}$ is a normed algebra with the usual quotient norm. In this section, we show that $\mathcal{A}[Z]/\mathcal{I}[Z]$ is a pre-C$^\ast$-algebra.  Also using Five Lemma and Theorem $\ref{theorem02}$ below, we will show that ${\mathcal{A}[Z]}/{\mathcal{I}[Z]}$ is isometric $\ast$-isomorphic to $({\mathcal{A}}/{\mathcal{I}})[Z]$.  This implies that the functor $\mathfrak{P}$ is, in fact, normal exact.

We remind that the Five Lemma in homological algebra (see, e.g., \cite{8}) says that  in the commutative diagram
\begin{center}
\begin{tikzpicture}[
scale=0.5,
ca/.style={shorten <=28.3pt},
cb/.style={shorten <=2pt},]
\matrix (m) [matrix of math nodes,row sep=4em,column sep=4em,minimum width=2em]
 {
A_1&A_2&A_3&A_4&A_5 \\
B_1&B_2&B_3&B_4&B_5 \\ };
   \path[-stealth]

   (m-1-1) edge node [right] {$t_1$} (m-2-1)
   (m-1-2) edge node [right] {$t_2$} (m-2-2)
   (m-1-3) edge node [right] {$t_3$} (m-2-3)
   (m-1-4) edge node [right] {$t_4$} (m-2-4)
   (m-1-5) edge node [right] {$t_5$} (m-2-5)
   (m-1-1) edge (m-1-2)
   (m-1-2) edge (m-1-3)
   (m-1-3) edge (m-1-4)
   (m-1-4) edge (m-1-5)
   (m-2-1) edge (m-2-2)
   (m-2-2) edge (m-2-3)
   (m-2-3) edge (m-2-4)
   (m-2-4) edge (m-2-5)
   ;
\end{tikzpicture}
\end{center}
of commutative $R$-modules with exact rows  if $t_1, t_2, t_4$ and $t_5$ are isomorphisms, so is  $t_3$.

\begin{definition}\label{definition03}\rm{
\cite{2} The exact sequence
$$\cdots\longrightarrow A_n\overset{f_n}{\longrightarrow}A_{n+1}\overset{f_{n+1}}{\longrightarrow}A_{n+2}\longrightarrow\cdots$$
of normed spaces with contraction $f_n$ ($\|f_n\|\leq 1$ for any $n$) is called normal exact if the induced map $A_n/\ker f_n\longrightarrow f_n(A_n)$ defined by $x+\ker f_n \longmapsto f_n(x)$, is an isometry. Note that any short exact sequence of C$^\ast$-algebras is normal exact.
}
\end{definition}

The following theorem is the main one   in {\cite{2}}.
\begin{theorem}  \label{theorem02}
 Suppose that
$$0\longrightarrow Y\overset{i}{\longrightarrow}X\overset{p}{\longrightarrow}Z\longrightarrow 0$$
is a normal exact sequence of normed spaces. Then
$$0\longrightarrow\hat{Y}\overset{\hat{i}}{\longrightarrow}\hat{X}\overset{\hat{p}}{\longrightarrow}\hat{Z}\longrightarrow 0 $$
is a normal exact sequence of corresponding completion Banach spaces.
\end{theorem}

\begin{theorem}\label{theorem03}
Let $\mathcal{I}$ be a closed ideal of a C$^\ast$-algebra $\mathcal{A}$. Then $\mathcal{A}[Z]/\mathcal{I}[Z]$ is a pre-C$^\ast$-algebra.

\noindent \begin{proof}
We first show that
\begin{description}
\item[(i)]If $(u_\lambda)_{\lambda\in\Lambda}$ is an approximate unit for $\mathcal{A}$, then $(u_\lambda)_{\lambda\in\Lambda}$ is also an approximate unit for $\mathcal{A}[Z]$;
\item[(ii)]If $(u_\lambda)_{\lambda\in\Lambda}$ is an approximate unit for $\mathcal{I}$, then for any $F(Z)\in\mathcal{A}[Z]$ we have
\begin{center}
\begin{tabular}{lll}
$\|F(Z)+\mathcal{I}[Z]\|$&$=$&$\lim_\lambda\|F(Z)-u_\lambda F(Z)\|_K$\\[0.2cm]
&$=$&$\lim_\lambda\|F(Z)-F(Z)u_\lambda\|_K$.\\[0.2cm]
\end{tabular}
\end{center}
\end{description}
To prove (i), let $F(Z)=\sum_{n=0}^{\infty}a_nZ^n\in\mathcal{A}[Z]$ and $\varepsilon >0$ be given. Since $\sum_{n=0}^{\infty}\|a_n\|<\infty$, there is a positive integer $N$ such  that $\sum_{n=N+1}^{\infty}2\|a_n\|<\varepsilon$. Now for any $\lambda\in\Lambda$ we have

\begin{center}
\begin{tabular}{lll}
$\|F(Z)-u_\lambda F(Z)\|_K$&$=$&$\|\sum_{n=0}^{\infty}(a_n-u_\lambda a_n)Z^n\|_K$\\[0.2cm]
&$\leq$&$\sum_{n=0}^{\infty}\|a_n-u_\lambda a_n\|$\\[0.2cm]
&$=$&$\sum_{n=0}^{N}\|a_n-u_\lambda a_n\|+\sum_{n=N+1}^{\infty}\|a_n-u_\lambda a_n\|$\\[0.2cm]
&$<$&$\sum_{n=0}^{N}\|a_n-u_\lambda a_n\|+\varepsilon$.\\[0.2cm]
\end{tabular}
\end{center}
Therefore
$$\lim_\lambda\sup \|F(Z)-u_\lambda F(Z)\|_K \leq \varepsilon.$$
Since $\varepsilon >0$ was arbitrary, we have
$$\lim_\lambda\|F(Z)-u_\lambda F(Z)\|_K =0.$$
Similarly we get
$$\lim_\lambda\|F(Z)- F(Z)u_\lambda\|_K =0.$$
To prove (ii) let
$$\alpha =\|F(Z)+\mathcal{I}[Z]\|=\inf\{\|F(Z)+H(Z)\|_K: H(Z)\in\mathcal{I}[Z]\}.$$
Let $\varepsilon >0$ be given. There exists an element $G(Z)\in\mathcal{I}[Z]$ such that $\|F(Z)-G(Z)\|_K<\alpha+\varepsilon$. We have
\begin{center}
\begin{tabular}{lll}
$\alpha$&$\leq$&$\|F(Z)-F(Z)u_\lambda\|_K$\\[0.2cm]
&$\leq$&$\|(F(Z)-G(Z))-(F(Z)-G(Z))u_\lambda\|_K+\|G(Z)-G(Z)u_\lambda\|_K$\\[0.2cm]
&$=$&$\|(F(Z)-G(Z))(1-u_\lambda)\|_K+\|G(Z)-G(Z)u_\lambda\|_K$\\[0.2cm]
&$\leq$&$\|F(Z)-G(Z)\|_K+\|G(Z)-G(Z)u_\lambda\|_K$\\[0.2cm]
&$<$&$\alpha +\varepsilon +\|G(Z)-G(Z)u_\lambda\|_K$.\\[0.2cm]
\end{tabular}
\end{center}
Now by Part (i) we have
\begin{center}
\begin{tabular}{lllll}
$\alpha$&$\leq$&$\lim_\lambda\inf\|F(Z)-F(Z)u_\lambda\|_K$&$\leq$&$\alpha +\varepsilon$,\\[0.2cm]
$\alpha$&$\leq$&$\lim_\lambda\sup\|F(Z)-F(Z)u_\lambda\|_K$&$\leq$&$\alpha +\varepsilon$.\\[0.2cm]
\end{tabular}
\end{center}
Since $\varepsilon >0$ was arbitrary, we have $\alpha =\lim_\lambda\|F(Z)-F(Z)u_\lambda\|_K$. Similarly, $\alpha =\lim_\lambda\|F(Z)-u_\lambda F(Z)\|_K.$

To prove the theorem let $(u_\lambda)_{\lambda\in\Lambda}$ be an approximate unit for $\mathcal{I}$. If $F(Z)\in\mathcal{A}[Z]$ and $G(Z)\in\mathcal{I}[Z]$, by Parts (i), (ii) and Proposition $\ref{proposition05}$ (i) we have
\begin{center}
\begin{tabular}{lll}
$\|F(Z)+\mathcal{I}[Z]\|^2$&$=$&$\lim_\lambda\|F(Z)-F(Z)u_\lambda\|_K^2$\\[0.2cm]
&$=$&$\lim_\lambda\|(1-u_\lambda)F^\ast(Z)F(Z)(1-u_\lambda)\|_K$\\[0.2cm]
&$\leq$&$\lim_\lambda\|(1-u_\lambda)(F^\ast(Z)F(Z)+G(Z))(1-u_\lambda)\|_K$\\[0.2cm]
&$+$&$\lim_\lambda\|(1-u_\lambda)G(Z)(1-u_\lambda)\|_K$\\[0.2cm]
&$\leq$&$\|F^\ast(Z)F(Z)+G(Z)\|_K$.\\[0.2cm]
\end{tabular}
\end{center}
Therefore
$$\|F(Z)+\mathcal{I}[Z]\|^2\leq\|F^\ast(Z)F(Z)+\mathcal{I}[Z]\|$$
and consequently we  get the equality
$$\|F(Z)+\mathcal{I}[Z]\|^2=\|F^\ast(Z)F(Z)+\mathcal{I}[Z]\|,$$
which completes the proof.
\end{proof}
\end{theorem}

Now we are ready to show that the functor $\mathfrak{P}$ is normal exact.

\begin{theorem}\label{theorem04}
The functor $\mathfrak{P}$ is normal exact.

\noindent \begin{proof}
 Let $\mathcal{I}$ be a closed ideal of a C$^\ast$-algebra $\mathcal{A}$. First we show that there exists an isometric $\ast$-isomorphism between $\mathcal{A}[Z]/\mathcal{I}[Z]$ and $(\mathcal{A}/\mathcal{I})[Z]$.
 Define $T:\mathcal{A}[Z]/\mathcal{I}[Z]\longrightarrow(\mathcal{A}/\mathcal{I})[Z]$ by
$$ T(\sum_{n=0}^{\infty}a_nZ^n+\mathcal{I}[Z])=\sum_{n=0}^{\infty}(a_n+\mathcal{I})Z^n, $$
for all $\sum_{n=0}^{\infty}a_nZ^n\in\mathcal{A}[Z]$.
 It is clear that $T$ is well defined, linear and preserves the involution. We are going to show  that
\textbf{(a)} $T$ is injective, \textbf{(b)} $T$ is surjective, \textbf{(c)} $T$ is a contraction, and \textbf{(d)} $T$ is an isometry.
We proceed as follows:
\begin{description}
\item[(a)] If $F(Z)=\sum_{n=0}^{\infty}a_nZ^n\in\mathcal{A}[Z]$ with $T(F)=\mathcal{I}$, then $\sum_{n=0}^{\infty}(a_n+\mathcal{I})Z^n=\mathcal{I}$, i.e., $a_n\in\mathcal{I}$ for $n=0, 1, 2,\cdots$. Therefore $F(Z)\in\mathcal{I}[Z]$ and so $T$ is injective.
\item[(b)] Let $G=\sum_{n=0}^{\infty}(a_n+\mathcal{I})Z^n\in(\mathcal{A}/\mathcal{I})[Z]$. For each $n=0, 1, 2,\cdots$ there is an element $b_n\in\mathcal{I}$ such that $\|a_n+b_n\|<\|a_n+\mathcal{I}\|+2^{-n}$. Let $\sum_{n=0}^{\infty}c_nZ^n$, where $c_n=a_n+b_n$ for each $n=0, 1, 2,\cdots$. Since $\sum_{n=0}^{\infty}\|a_n+\mathcal{I}\|<\infty$ we have $F(Z)\in\mathcal{A}[Z]$. Therefore
\begin{center}
\begin{tabular}{lll}
$T(F(Z)+\mathcal{I}[Z])$&$=$&$\sum_{n=0}^{\infty}(c_n+\mathcal{I})Z^n$\\[0.2cm]
&$=$&$\sum_{n=0}^{\infty}(a_n+\mathcal{I})Z^n$\\[0.2cm]
&$=$&$G$,\\[0.2cm]
\end{tabular}
\end{center}
that is $T$ is surjective.
\item[(c)] Let $F(Z)=\sum_{n=0}^{\infty}a_nZ^n\in\mathcal{A}[Z]$.  Then

\begin{center}
\begin{tabular}{lll}
$\|T(F(Z)+\mathcal{I}[Z])\|$&$=$&$\|\sum_{n=0}^{\infty}(a_n+\mathcal{I})Z^n\|_K$\\[0.2cm]
&$=$&$\sup_{t\in K}\|\sum_{n=0}^{\infty}(a_n+\mathcal{I})t^n\|$\\[0.2cm]
&$=$&$\sup_{t\in K}\|\sum_{n=0}^{\infty}a_nt^n+\mathcal{I}\|$\\[0.2cm]
&$=$&$\sup_{t\in K}\inf_{b\in\mathcal{I}}\|\sum_{n=0}^{\infty}a_nt^n+b\|$\\[0.2cm]
&$\leq$&$\inf_{b\in\mathcal{I}}\sup_{t\in K}\|\sum_{n=0}^{\infty}a_nt^n+b\|$\\[0.2cm]
&$=$&$\inf_{G(Z)\in\mathcal{I}[Z]}\sup_{t\in K}\|\sum_{n=0}^{\infty}a_nt^n+G(t)\|$\\[0.2cm]
&$=$&$\inf_{G}\|F(Z)+G(Z)\|_K$\\[0.2cm]
&$=$&$\|F(Z)+\mathcal{I}[Z]\|$,\\[0.2cm]
\end{tabular}
\end{center}
that is $T$ is a contraction. (Note that $\sup\inf f\leq \inf\sup f$ for every real valued function $f$ in two variables)
\item[(d)] Suppose that  $\widehat{(\mathcal{A}[Z]/\mathcal{I}[Z])}$ is the completion of $\mathcal{A}[Z]/\mathcal{I}[Z]$ with respect to the  quotient norm and
$$\hat{T}:(\widehat{\mathcal{A}[Z]/\mathcal{I}[Z]})\longrightarrow [\mathcal{A}/\mathcal{I}]_K,$$
is the extension of $T$. By Theorem $\ref{theorem03}$, $\hat{T}$ is a $\ast$-homomorphism of C$^\ast$-algebras. Now we show that $\hat{T}$ is a $\ast$-isomorphism. The diagram
\begin{center}
\begin{tikzpicture}[
scale=0.5,
ca/.style={shorten <=28.3pt},
cb/.style={shorten <=2pt},]
\matrix (m) [matrix of math nodes,row sep=3em,column sep=3em,minimum width=2em]
 {
\sum_{n=0}^{\infty}a_nZ^n&\sum_{n=0}^{\infty}a_nZ^n+\mathcal{I}[Z]\\
\sum_{n=0}^{\infty}a_nZ^n&\sum_{n=0}^{\infty}(a_n+\mathcal{I})Z^n \\ };
   \path[-stealth]
   (m-1-1) edge [|->]  (m-2-1)
   (m-1-2) edge [|->] node [right] {$T$}(m-2-2)
   (m-2-1) edge [|->] node [above] {$\tilde{p}$}(m-2-2)
   (m-1-1) edge [|->] node [above] {$p'$}(m-1-2)
   ;
\end{tikzpicture}
\end{center}
shows that the diagram
\begin{center}
\begin{tikzpicture}[
scale=0.5,
ca/.style={shorten <=28.3pt},
cb/.style={shorten <=2pt},]
\matrix (m) [matrix of math nodes,row sep=3em,column sep=3em,minimum width=2em]
 {
0&\mathcal{I}[Z]&\mathcal{A}[Z]&\mathcal{A}[Z]/\mathcal{I}[Z]&0 \\
0&\mathcal{I}[Z]&\mathcal{A}[Z]&(\mathcal{A}/\mathcal{I})[Z]&0 \\ };
   \path[-stealth]
   (m-1-2) edge [-,double] (m-2-2)
   (m-1-3) edge [-,double] (m-2-3)
   (m-1-4) edge node [right] {$T$} (m-2-4)
   (m-1-1.east|-m-1-2) edge (m-1-2)
   (m-1-2) edge [cb] (m-1-3)

   (m-1-3) edge [-right hook,ca] (m-1-2)
   (m-2-3) edge [-right hook,ca] (m-2-2)
   (m-1-3) edge node [above] {$p'$} (m-1-4)
   (m-1-4.east|-m-1-5) edge (m-1-5)
   (m-2-1.east|-m-2-2) edge (m-2-2)
   (m-2-2) edge [cb] (m-2-3)
   (m-2-3) edge node [above] {$\tilde{p}$} (m-2-4)
   (m-2-4.east|-m-2-5) edge (m-2-5)
   ;
\end{tikzpicture}
\end{center}
of pre-C$^\ast$-algebras is commutative,
where $p'$ is the quotient map and $\tilde{p}$ is the map induced by the projection $\mathcal{A}\overset{p}{\longrightarrow}\mathcal{A}/\mathcal{I}$ (see Definition $\ref{definition02}$). The exactness of two rows follow from Proposition $\ref{proposition08}$ (vii). Now the commutative diagram
\begin{center}
\begin{tikzpicture}[
scale=0.5,
ca/.style={shorten <=28.3pt},
cb/.style={shorten <=2pt},]
\matrix (m) [matrix of math nodes,row sep=3em,column sep=3em,minimum width=2em]
 {
0&[][\mathcal{I}]_K&[][\mathcal{A}]_K&(\widehat{\mathcal{A}[Z]/\mathcal{I}[Z]})&0 \\
0&[][\mathcal{I}]_K&[][\mathcal{A}]_K&[][\mathcal{A}/\mathcal{I}]_K&0 \\ };
   \path[-stealth]
   (m-1-2) edge [-,double] (m-2-2)
   (m-1-3) edge [-,double] (m-2-3)
   (m-1-4) edge node [right] {$\hat{T}$} (m-2-4)
   (m-1-1.east|-m-1-2) edge (m-1-2)
   (m-1-2) edge [cb] (m-1-3)
   (m-1-3.east|-m-1-3) edge node [above] {$\hat{p'}$} (m-1-4.west|-m-1-3)
   (m-1-4.east|-m-1-5) edge (m-1-5)
   (m-2-1.east|-m-2-2) edge (m-2-2)
   (m-2-2.east|-m-2-3) edge [cb] (m-2-3)
   (m-2-3.east|-m-2-3) edge node [above] {$\hat{p}$} (m-2-4.west|-m-2-4)
   (m-2-4.east|-m-2-5) edge (m-2-5)
      (m-1-3.east|-m-1-2) edge [-right hook,ca] (m-1-2)
   (m-2-3.east|-m-2-2) edge [-right hook,ca] (m-2-2)
   ;
\end{tikzpicture}
\end{center}
of C$^\ast$-algebras have exact rows. In fact, the  exactness of first row is a consequence of Theorem $\ref{theorem02}$ and the second one follows from Theorem $\ref{theorem01}$. Applying Five Lemma for commutative diagram
\begin{center}
\begin{tikzpicture}[
scale=0.5,
ca/.style={shorten <=28.3pt},
cb/.style={shorten <=2pt},
cc/.style={shorten >=90pt},]
\matrix (m) [matrix of math nodes,row sep=3em,column sep=3em,minimum width=2em]
 {
[\mathcal{I}]_K&[][\mathcal{A}]_K&(\widehat{\mathcal{A}[Z]/\mathcal{I}[Z]})&0&0\\
[][\mathcal{I}]_K&[][\mathcal{A}]_K&[][\mathcal{A}/\mathcal{I}]_K&0&0 \\ };
   \path[-stealth]
   (m-1-1) edge [-,double] node [left] {$t_1$}(m-2-1)
   (m-1-2) edge [-,double] node [left] {$t_2$}(m-2-2)
   (m-1-3) edge node [left] {$t_3=\hat{T}$}(m-2-3)
   (m-1-4) edge node [left] {$t_4$} (m-2-4)
   (m-1-5) edge node [left] {$t_5$} (m-2-5)
   (m-1-1) edge [cb] (m-1-2)
   (m-1-2) edge node [above] {$\hat{p'}$} (m-1-3.west|-m-1-2)
   (m-1-3.east|-m-1-4) edge (m-1-4.west|-m-1-4)
   (m-1-4) edge (m-1-5)
   (m-2-1) edge [cb] (m-2-2)
   (m-2-2) edge node [above] {$\hat{p}$}  (m-2-3)
   (m-2-3.east|-m-2-4) edge (m-2-4.west|-m-2-4)
   (m-2-4) edge (m-2-5)
      (m-1-2) edge [-right hook,ca] (m-1-1)
   (m-2-2) edge [-right hook,ca] (m-2-1)
   ;
\end{tikzpicture}
\end{center}
with exact rows shows that $\hat{T}$ is a $\ast$-isomorphism. This implies, particularly, that $T$ is an isometry. Now consider the short exact sequence of C$^\ast$-algebras
$$0\longrightarrow\mathcal{I}\overset{i}{\hookrightarrow}\mathcal{A}\overset{g}{\longrightarrow}\mathcal{B}\longrightarrow 0.$$
Applying functor $\mathfrak{P}$ we get a short exact sequence of pre-C$^\ast$-algebras
\begin{equation}\label{equation13}
0\longrightarrow\mathcal{I}[Z]\overset{\tilde{i}}{\hookrightarrow}\mathcal{A}[Z]\overset{\tilde{g}}{\longrightarrow}\mathcal{B}[Z]\longrightarrow 0.
\end{equation}
Note that we have the $\ast$-isomorphism $g_1:\mathcal{A}/\mathcal{I}\longrightarrow\mathcal{B}$, induced by $g$. By Part (d) we have the composition of isometric $\ast$-isomorphism of pre-C$^\ast$-algebras
$$\mathcal{A}[Z]/\mathcal{I}[Z]\overset{T}{\longrightarrow}(\mathcal{A}/\mathcal{I})(Z)\overset{\tilde{g}_1}{\longrightarrow}\mathcal{B}[Z]$$
such that $$\sum_{n=0}^{\infty}a_nZ^n+\mathcal{I}[Z]\mapsto\sum_{n=0}^{\infty}(a_n+\mathcal{I})Z^n\mapsto\sum_{n=0}^{\infty}g(a_n)Z^n.$$ 
That, is the induced map $\mathcal{A}[Z]/\mathcal{I}[Z]\longrightarrow\mathcal{B}[Z]$ by $\tilde{g}$ is an isometry. Therefore $(\ref{equation13})$ is a normal exact sequence of pre-C$^\ast$-algebras.

\end{description}
\end{proof}
\end{theorem}

From (c) and (d) of Theorem $\ref{theorem04}$ we have the following.
\begin{corollary}\label{corollary04}
Suppose that $\mathcal{I}$ is a closed ideal of a C$^\ast$-algebra $\mathcal{A}$ and $a_0, a_1,\\ a_2,\cdots$, is a sequence in $\mathcal{A}$ such that $\sum_{n=0}^{\infty}\|a_n\|<\infty$. Then
$$\inf_{b\in\mathcal{I}}\sup_{t\in K}\|\sum_{n=0}^{\infty}a_nt^n+b\|=\sup_{t\in K}\inf_{b\in\mathcal{I}}\|\sum_{n=0}^{\infty}a_nt^n+b\|.$$
\end{corollary}

\section{Cauchy extension $[\mathcal{A}]_K$ as C$^*$-subalgebra of \\ $C_b(K,\mathcal{A})$ }\label{section06}
In this section, we characterize the Cauchy extensions of C$^\ast$-algebras as  C$^\ast$-valued function spaces. Using the obtained characterization, we give some  results on the  Cauchy extensions of C$^\ast$-algebras.

Recall that for a C$^\ast$-algebra $\mathcal{A}$ and a topological space $X$, $C_b(X,\mathcal{A})$ is the set of all bounded continuous functions from $X$ to $\mathcal{A}$. The addition, scalar multiplication and the product on $C_b(X,\mathcal{A})$ are defined pointwise. The involution can be defined as  $\alpha ^\ast (x)=(\alpha (x)) ^\ast$, for all $\alpha\in C_b(X,\mathcal{A})$ and $x\in X$. Furthermore, defining  $\|\alpha\|_\infty=\sup_{x\in X}\|\alpha (x)\|$ for all $\alpha\in C_b(X,\mathcal{A})$, the algebra $C_b(X,\mathcal{A})$ becomes  a C$^\ast$-algebra. If $X$ is a locally compact Hausdroff space, then $C_0(X,\mathcal{A})$ consisting of  all continuous functions $f\in C_b(X,\mathcal{A})$ vanishing at infinity is a C$^\ast$-subalgebra of $C_b(X,\mathcal{A})$  (see \cite[p.37]{9} ). If $X$ is a compact Hausdorff space, then $C_b(X,\mathcal{A})=C_0(X,\mathcal{A})=C(X,\mathcal{A}).$

It is easy to see that for C$^\ast$-algebras $\mathcal{A}_1, \mathcal{A}_2,\cdots ,\mathcal{A}_n$, we have
\begin{equation}\label{equation14}
C_b(X,\mathcal{A}_1\oplus\cdots\oplus\mathcal{A}_n)\cong C_b(X,\mathcal{A}_1)\oplus\cdots\oplus C_b(X,\mathcal{A}_n).
\end{equation}
In particular, if $\mathcal{A}=\mathbb{C}$, we use $C(X), C_b(X)$ and $C_0(X)$ for $C(X, \mathbb{C}), C_b(X, \mathbb{C})$ and $C_0(X, \mathbb{C})$, respectively.
Recall that a C$^\ast$-algebra $\mathcal{A}$ is called nuclear if for each C$^\ast$-algebra $\mathcal{B}$, there is a unique C$^\ast$-norm on tensor product $\mathcal{A}\otimes\mathcal{B}$. An ideal $\mathcal{I}$ of a C$^\ast$-algebra $\mathcal{A}$ is called essential if  $a\mathcal{I}=0$ implies that  $a=0$.
\begin{theorem}\label{theorem05}
Suppose that $\mathcal{A}$ and $\mathcal{B}$ are two C$^\ast$-algebras and $K\subseteq  J=[-1,1]$ such that $0$ is a limit point of $K$. Then
\begin{description}
\item[(i)] $[\mathcal{A}]_K$ is $\ast$-isomorphic to a C$^\ast$-subalgebra of $C_b(K,\mathcal{A})$;
\item[(ii)] If $K$ is a compact interval, then $[\mathcal{A}]_K\cong C(K,\mathcal{A})$;
\item[(iii)] $[\mathcal{A}]_K\cong\{f|_{K}:f\in C([-1,1],\mathcal{A})\}$;
\item[(iv)] If $K$ is compact, then $[\mathcal{A}]_K\cong C(K,\mathcal{A})$. Furthermore, $[\mathcal{A}\otimes\mathcal{B}]_K\cong[\mathcal{A}]_K\otimes\mathcal{B}\cong \mathcal{A}\otimes[\mathcal{B}]_K$;
\item[(v)] There is a closed ideal $\mathcal{I}_K$ of $[\mathcal{A}]_J$ such that $[\mathcal{A}]_J/\mathcal{I}_K\cong [\mathcal{A}]_K$;
\item[(vi)] $\mathcal{A}$ is nuclear if and only if $[\mathcal{A}]_K$ is nuclear;
\item[(vii)] $\mathcal{I}$ is an essential ideal of $\mathcal{A}$ if and only if $[\mathcal{I}]_K$ is an essential ideal of $[\mathcal{A}]_K$;
\item[(viii)] If $0\notin K$ and $K$ is a locally compact subspace of $J$ such that $K'=K\cup\{0\}$ is compact, then $[\mathcal{A}]_K\cong C(K',\mathcal{A})$. If $\mathcal{A}$ is finite dimensional, then $M(\hat{\mathcal{A}}_1)\cong C_b(K,\mathcal{A})$, where $M(\hat{\mathcal{A}}_1)$ is the multiplier algebra of $\hat{\mathcal{A}}_1$;
\item[(ix)] $\mathcal{A}\cong\mathcal{B}$ if and only if $[\mathcal{A}]_K\cong[\mathcal{B}]_K$ for any compact set $K$.
\end{description}
\noindent \begin{proof}

\item[(i)] It is clear that for any sequence $(a_n)$ in $\mathcal{A}$ with $\sum_{n=0}^{\infty}\|a_n\|<\infty$ the summation $f(t)=\sum_{n=0}^{\infty}a_nt^n$, where $t\in K$ defines a function from $K$ to $\mathcal{A}$. Denote the set of all such functions by $\mathcal{A}(K)$. It is clear that $f$ is a bounded continuous function on $K$ and $\mathcal{A}(K)$ is a $\ast$-subalgebra of $C_b(K,\mathcal{A})$. Now the map $T:\mathcal{A}(K)\longrightarrow\mathcal{A}[Z]$ defined by $T(\sum_{n=0}^{\infty}a_nt^n)=\sum_{n=0}^{\infty}a_nZ^n$ is an isometric  $\ast$-isomorphism. That is $[\mathcal{A}]_K$ is $\ast$-isomorphic to a C$^\ast$-subalgebra of $C_b(K,\mathcal{A})$.
\item[(ii)] For the case that $\mathcal{A}=\mathbb{C}$, since $\mathbb{C}(K)$ is a self-adjoint subalgebra of $C(K)$ which separate points of $K$ and contains the constant functions one can see, by Stone-Weierstrass Theorem (see \cite[p.165]{4}), that $[\mathbb{C}]_K\cong C(K)$. Now for any C$^\ast$-algebra $\mathcal{A}$ and any compact interval $K$ one can use approximate Berstein Theorem (see \cite[p.182]{12} ), as follows: We may assume  that $K=[0,1]$. Let $f\in C(K,\mathcal{A})$. Because $f$ is uniformly continuous (see \cite[p.60]{11} ),  define the Bernstein Polynomials  $$\beta_n(t)=\sum_{m=0}^{n}f(m/n)\binom{n}{m}t^m(1-t)^{n-m},$$ for any $t\in K$ and integer $n>0$. Note that $\beta_n\in\mathcal{A}(K)$ for any $n=1, 2, 3,\cdots$. By a  similar argument as  in the proof of the Berstein Theorem, we see that $\beta_n$ is convergent uniformly to $f$. This shows that $[\mathcal{A}]_K\cong C(K,\mathcal{A})$.
\item[(iii)] Define $T:\mathcal{A}(J)\longrightarrow\mathcal{A}(K)$ by $T(f)=f|_{K}$, for each $f\in\mathcal{A}(J)$. It is clear that $T$ is a bijective bounded linear operator. We claim that the extension $\hat{T}:[\mathcal{A}]_J\longrightarrow [\mathcal{A}]_K$ is surjective. Note that Parts (i) and (ii) imply that $\hat{T}$ is of the form $\hat{T}(g)=g|_{K}$ for all $g\in[\mathcal{A}]_J$. Suppose that $H, G:[\mathcal{A}]_K\longrightarrow\mathcal{B}$ are two $\ast$-homomorphisms such that $G\circ\hat{T}=H\circ\hat{T}$. This implies that $H\circ\hat{T}|_{\mathcal{A}(J)}=G\circ\hat{T}|_{\mathcal{A}(J)}$ or $H|_{\mathcal{A}(K)}=G|_{\mathcal{A}(K)}$. Since $\mathcal{A}(K)\cong\mathcal{A}[Z]$ is dense in $[\mathcal{A}]_K$, then $H=G$. Hence  $\hat{T}$ is surjective (see \cite{6}). By (ii) we have $[\mathcal{A}]_J\cong C(J,\mathcal{A})$ and therefore $[\mathcal{A}]_K\cong\{f|_K:f\in C(J,\mathcal{A})\}$.
\item[(iv)] By Tietze's Theorem (\cite[Theorem 4.1]{17}), any continuous function $f:K\longrightarrow \mathcal{A}$ has a continuous extension $f_1:J\longrightarrow\mathcal{A}$. This fact together with Part (iii) show that $[\mathcal{A}]_K\cong C(K,\mathcal{A})$. From (\cite[II.6.4.4]{13} ) we have $C(K,\mathcal{A})\cong C(K)\otimes\mathcal{A}$ and therefore
 $$[\mathcal{A}\otimes\mathcal{B}]_K\cong C(K)\otimes (\mathcal{A}\otimes\mathcal{B})\cong [\mathcal{A}]_K\otimes \mathcal{B}\cong\mathcal{A}\otimes [\mathcal{B}]_K.$$
\item[(v)] Let $\hat{T}:[\mathcal{A}]_J\longrightarrow [\mathcal{A}]_K$ be the given surjective $\ast$-homomorphism in Part (iii). If $\mathcal{I}_K=\ker\hat{T}$, then $[\mathcal{A}]_J/\mathcal{I}_K\cong [\mathcal{A}]_K$. In fact, $[\mathcal{A}]_J$ is an extension of any Cauchy extension $[\mathcal{A}]_K$.
\item[(vi)] Let $\mathcal{A}$ be nuclear. By Part (ii) we have $[\mathcal{A}]_J\cong C(J,\mathcal{A})$. Since $C(J)$ is nuclear (see \cite[Theorem 6.4.15]{9}) and $C(J,\mathcal{A})\cong C(J)\otimes\mathcal{A}$ (see \cite[II.6.4.4]{13}) we imply that $[\mathcal{A}]_J$ is nuclear (see \cite[IV.3.1.1]{13}). Since every closed ideal of a nuclear C$^\ast$-algebra is nuclear (see \cite[II.9.6.3]{13}), then $\hat{\mathcal{A}}_1$ is nuclear. In particular, the closed ideal $\mathcal{I}_K$ (given in part (v)) is nuclear. Since $[\mathcal{A}]_J/\mathcal{I}_K$ is nuclear (see \cite[IV 3.1.13]{13}), Part (v) implies that $[\mathcal{A}]_K$ is also nuclear. Conversely, if $[\mathcal{A}]_K$ is nuclear, then the ideal $\hat{\mathcal{A}}_1$ is nuclear. By $(\ref{equation08})$ we have $\mathcal{A}\cong [\mathcal{A}]_K/\hat{\mathcal{A}}_1$ which shows that $\mathcal{A}$ is  nuclear too.
\item[(vii)] Let $\mathcal{I}$ be an essential ideal of $\mathcal{A}$. By Part (i) we can consider $[\mathcal{A}]_K$ as a C$^\ast$-subalgebra of $C_b(K,\mathcal{A})$. Choose $G:K\longrightarrow\mathcal{A}$ in $[\mathcal{A}]_K$ such that $fG=Gf=0$ for any $f:K\longrightarrow\mathcal{I}\in [\mathcal{I}]_K$. For any $t\in K$ we have $f(t)G(t)=G(t)f(t)=0$. Let $b$ be an arbitrary element in $\mathcal{I}$ and let $f_b:K\longrightarrow\mathcal{I}$ be a constant function with value $f_b(t)=b$. Now for any $t\in K$ we have
$$f_b(t)G(t)=G(t)f_b(t)=0$$
or
$$bG(t)=G(t)b=0.$$
This implies that $G(t)=0$ for all $t\in K$. Therefore $[\mathcal{I}]_K$ is an essential ideal of $[\mathcal{A}]_K$. The converse statement can be proved similarly.
\item[(viii)] Suppose that $\mathbb{C}_1(K)=\{f\in\mathbb{C}(K):f(0)=0\}$, where $\mathbb{C}(K)$ is as given in Part (i). For $f\in \mathbb{C}_1(K)$ and $\varepsilon >0$, suppose that $X=\{t\in K:|f(t)|\geq \varepsilon\}$ and $x$ is a limit point of $X$. Then $x\neq 0$ and $x$ is a limit point of $K'$, and  therefore $x\in K$. This implies that $X$ is compact. That is $f$ vanishes at infinity, so $\mathbb{C}_1(K)\subseteq  C_0(K)$.  Now suppose that $0\neq a\in\mathbb{C}$ and $g(x)=xa$ for all $x\in K$. Then $g\in\mathbb{C}_1(K)$ and for any $t\in K$ we have $g(t)\neq 0$. In addition, for any $t_1\neq t_2$ in $K$, $g(t_1)\neq g(t_2)$, that is, $\mathbb{C}_1(K)$ strongly separates points of $K$. It is clear that $\mathbb{C}_1(K)$ is self-adjoint. By the Stone-Weierstrass Theorem (see \cite[p.151]{10}) we have $\hat{\mathbb{C}}_1\cong C_0(K)$ and therefore $[\mathbb{C}]_K\cong \mathbb{C}\oplus C_0(K)\cong C(K')$ (see \cite[p.53]{13}). Parts (iii) and (iv) and the fact that $\|f\|_K=\|f\|_{K'}$ for any $f\in C(J,\mathcal{A})$ imply  that the map $f|_K\mapsto f|_{K'}$ is a $\ast$-isomorphism between $[\mathcal{A}]_K$ and $C(K',\mathcal{A})$. Now suppose that $\mathcal{A}$ is a finite dimensional C$^\ast$-algebra. By (\cite[p.194]{9}) we have
\begin{equation}\label{equation15}
\mathcal{A}\cong M _{n_{1}}(\mathbb{C})\oplus M_{n_{2}}(\mathbb{C})\oplus\cdots\oplus M_{n_{m}}(\mathbb{C}).
\end{equation}
We first  show that for any positive integer $n, \widehat{(M_n(\mathbb{C}))_1}\cong M_n(\hat{\mathbb{C}}_1)$. To see this, note that the completion of $\mathbb{C}_1(K)$ is $\ast$-isomorphic to $C_0(K)$. Now the map $G:(M_n(\mathbb{C}))_1(K)\longrightarrow M_n(\mathbb{C}_1(K))$ defined by $G(F)=(F_{ij})$, where
$$F(t)=\sum_{m=1}^{\infty}B_mt^m=(F_{ij}(t))$$
 and $F_{ij}\in\mathbb{C}_1(K)$, for any $i, j=1, 2,\cdots, n$ is an isometric $\ast$-isomorphism with  norm $\|(F_{ij})\|=\sup_{t\in K}\|F_{ij}(t)\|=\|F\|$. Suppose that $F=(F_{ij})\in M_n(C_0(K))$, then $F_{ij}\in C_0(K)$ for $i, j=1, 2,\cdots, n$. There exist sequences $(F_{mij})$ in $\mathbb{C}_1(K)$ for $i,j=1, 2,\cdots, n$ such that $F_{mij}\longrightarrow F_{ij}$ as $m\to\infty$ in norm $\|\cdot\|_K$. If $F:K\longrightarrow M_n(\mathbb{C})$ is a continuous function such that for any $t\in K, F(t)=(F_{ij}(t))$, then
\begin{center}
\begin{tabular}{lll}
$\|(F_{mij})-(F_{ij})\|$&$=$&$\sup_{t\in K}\|(F_{mij}(t))-F_{ij}(t)\|$\\[0.2cm]
&$\leq$&$\sup_{t\in K}\sum_{i, j}\|F_{mij}(t)-F_{ij}(t)\|$\\[0.2cm]
&$\leq$&$\sum_{i, j}\sup_{t\in K}\|F_{mij}(t)-F_{ij}(t)\|$.\\[0.2cm]
\end{tabular}
\end{center}
This implies that $(F_{mij})\longrightarrow (F_{ij})$ as $m\longrightarrow\infty$.
 Now by completion we see that
 \begin{equation}\label{equation16}
\widehat{(M_n(\mathbb{C}))}_1\cong M_n(C_0(K))\cong  M_n(\hat{\mathbb{C}}_1).
 \end{equation}
 Also we have clearly the $\ast$-isomorphism
 \begin{equation}\label{equation17}
M_n(C_b(K))\cong C_b(K,M_n(C)).
 \end{equation}
 Suppose that $\mathcal{B}, \mathcal{A}_1, \mathcal{A}_2,\cdots,\mathcal{A}_n$ are C$^\ast$-algebras. We have the following for the multipliers algebras (see \cite[p.84]{16})
\begin{equation}\label{equation18}
M(M_n(\mathcal{B}))\cong M_n(M(\mathcal{B}))
\end{equation}
\begin{equation}\label{equation19}
M(\mathcal{A}_1\oplus\mathcal{A}_2\oplus\cdots\oplus \mathcal{A}_n)\cong M(\mathcal{A}_1)\oplus M(\mathcal{A}_2)\oplus\cdots \oplus M(\mathcal{A}_n),
\end{equation}

(see \cite[p.155]{13}). We also have $M(C_0(K))\cong C_b(K)$ (see \cite[p.83]{9}). Now from $(\ref{equation14})-(\ref{equation19}),$ and Proposition \ref{proposition08} (ix), we have
$$\hat{\mathcal{A}}_1\cong M_{n_1}(C_0(K))\oplus M_{n_2}(C_0(K))\oplus\cdots\oplus M_{n_m}(C_0(K)).$$
  \begin{center}
\begin{tabular}{lll}
$M(\hat{\mathcal{A}}_1)$&$\cong$&$M(M_{n_1}(C_0(K)))\oplus\cdots\oplus M(M_{n_m}(C_0(K)))$\\[0.2cm]
&$\cong$&$C_b(K, M_{n_1}(\mathbb{C}))\oplus\cdots\oplus C_b(K, M_{n_m}(\mathbb{C}))$\\[0.2cm]
&$\cong$&$C_b(K, M_{n_1}(\mathbb{C})\oplus\cdots\oplus  M_{n_m}(\mathbb{C}))$\\[0.2cm]
&$\cong$&$C_b(K, \mathcal{A})$.\\[0.2cm]
\end{tabular}
\end{center}
\item[(ix)] If $\mathcal{A}\cong \mathcal{B}$, then $[\mathcal{A}]_K\cong[\mathcal{B}]_K$ by Proposition \ref{proposition08} (iv). Let $\varphi_n :[\mathcal{A}]_{K_n}\longrightarrow[\mathcal{B}]_{K_n}$ be a $\ast$-isomorphism between $[\mathcal{A}]_{K_n}$ and $[\mathcal{B}]_{K_n}$, where $K_n=[-1/n, 1/n]$
for $n=1,2,3,\cdots$. It is clear that $(K_n)$ is nested with $\bigcap_{n=1}^{\infty}K_n=\{0\}$. Now $([\mathcal{A}]_{K_n}, p_n)_{n=1}^{\infty}$ is a direct sequence of C$^\ast$-algebras, where each map 
$$p_n:[\mathcal{A}]_{K_{n}}\longrightarrow [\mathcal{A}]_{K_{n+1}}$$
 defined by $f|_{K_{n}}\mapsto f|_{K_{n+1}}$, for all $f\in[\mathcal{A}]_K$ is a 
$\ast$-homomorphism. Part (iv) and \cite[II.6.4.4]{13} show that
$$[\mathcal{A}]_{K_n}\cong C(K_n, \mathcal{A})\cong C(K_n)\otimes\mathcal{A},$$
for all $n$. Furthermore by \cite[II.9.6.5]{13} we have the direct limit
$$\lim_{\longrightarrow}[\mathcal{A}]_{K_n}\cong\lim_{\longrightarrow}(C(K_n)\otimes\mathcal{A})\cong(\lim_{\longrightarrow}C(K_n))\otimes\mathcal{A}\cong C(\{0\})\otimes\mathcal{A}\cong\mathbb{C}\otimes\mathcal{A}\cong\mathcal{A}.$$
From the commutative diagram
\begin{center}
\begin{tikzpicture}[
scale=0.5,
ca/.style={shorten <=28.3pt},
cb/.style={shorten <=2pt},]
\matrix (m) [matrix of math nodes,row sep=3em,column sep=3em,minimum width=2em]
 {
[\mathcal{A}]_{K_n}&[][\mathcal{B}]_{K_n}\\
[][\mathcal{A}]_{K_{n+1}}&[][\mathcal{B}]_{K_{n+1}} \\ };
   \path[-stealth]
   (m-1-1) edge [->] node [left] {$p_n$} (m-2-1)
   (m-1-2) edge [->] node [right] {$q_n$}(m-2-2)
   (m-2-1) edge [->] node [above]{$\varphi_{n+1}$}(m-2-2)
   (m-1-1) edge [->] node [above] {$\varphi_n$}(m-1-2)
   ;
\end{tikzpicture}
\end{center}
where $([\mathcal{B}]_{K_n},q_n)_{n=1}^{\infty}$ is the direct sequence defined by $q_n(\varphi_n(f))=\varphi_{n+1}(f|_{K_{n+1}})$, for any $f\in[\mathcal{A}]_{K_n}$, we conclude that
$$\mathcal{A}\cong \lim_{\longrightarrow}[\mathcal{A}]_{K_n}\cong \lim_{\longrightarrow}[\mathcal{B}]_{K_n}\cong\mathcal{B},$$
as desired.
\end{proof}
\end{theorem}

Any C$^\ast$-algebra of the form
$$\mathcal{B}=M_{n_1}(C[a_1,b_1])\oplus\cdots\oplus M_{n_m}(C[a_n,b_n])$$
where $a_i<b_i$ for $i=1, 2,\cdots,n$ are real numbers, is a Cauchy extension of some C$^\ast$-algebra. In fact
$$\mathcal{B}\cong M_{n_1}(C[-1,1])\oplus M_{n_2}(C[-1,1])\oplus\cdots\oplus M_{n_m}(C[-1,1]).$$
Therefore $\mathcal{B}\cong[\mathcal{A}]_J$, where $\mathcal{A}$ is the C$^\ast$-algebra defined in $(\ref{equation15})$.
\begin{corollary}\label{corollary05}
Suppose that $\mathcal{A}$ is a C$^\ast$-algebra and $\mathcal{I}$ is a closed ideal of $\mathcal{A}$. If $K=[0,1]$ and $F\in C(K,\mathcal{A})$, then
$$
\inf_{b\in \mathcal{I}}\sup_{t\in K}\|F(t)+b\|=\sup_{t\in K}\inf_{b\in \mathcal{I}}\|F(t)+b\|.
$$

\noindent \begin{proof}
Let $\varepsilon>0$ be  given. By Theorem $\ref{theorem05}$ (ii) there exists an element $F_n\in\mathcal{A}(K)$ such that $\sup_{t\in K}\|F(t)-F_n(t)\|<\varepsilon$. For any $t\in K$ we have
$$\|F(t)+b\|\leq\|F(t)-F_n(t)\|+\|F_n(t)+b\|<\varepsilon+\|F_n(t)+b\|.$$
On the other hand
$$\|F_n(t)+b\|\leq\|F_n(t)-F(t)\|+\|F(t)+b\|<\varepsilon +\|F(t)+b\|,$$ for any $t\in K$. By Corollary $\ref{corollary04}$ we have
 \begin{center}
\begin{tabular}{lll}
$\inf_{b\in \mathcal{I}}\sup_{t\in K}\|F(t)+b\|$&$\leq$&$\varepsilon +\inf_{b\in \mathcal{I}}\sup_{t\in K}\|F_n(t)+b\|$\\[0.2cm]
&$=$&$\varepsilon +\sup_{t\in K}\inf_{b\in \mathcal{I}}\|F_n(t)+b\|$\\[0.2cm]
&$\leq$&$2\varepsilon +\sup_{t\in K}\inf_{b\in \mathcal{I}}\|F(t)+b\|$.\\[0.2cm]
\end{tabular}
\end{center}
Since $\varepsilon>0$ was arbitrary, then
$$\inf_{b\in \mathcal{I}}\sup_{t\in K}\|F(t)+b\|\leq\sup_{t\in K}\inf_{b\in \mathcal{I}}\|F(t)+b\|.$$
This completes the proof.
\end{proof}
\end{corollary}


\begin{thebibliography}{10}
\bibitem{5}Arveson, W. Notes on extensions of $C\sp{\sp*}$ -algebras. Duke Math. J. \textbf{44} (1977), no. 2, 329--355.




\bibitem{12}Bartle, R.G., The elements of real analysis. Second edition. John Wiley $\&$ Sons, New York-London-Sydney, 1976.

\bibitem{13}Blackadar, B., Operator algebras. Theory of $C\sp *$ -algebras and Von Neumann algebras. Encyclopaedia of Mathematical Sciences, 122. Operator Algebras and Non-commutative Geometry, III. Springer-Verlag, Berlin, 2006.

\bibitem{16}Blecher, D.P.; Le Merdy, C. Operator algebras and their modules an operator space approach. London Mathematical Society Monographs. New Series, 30. Oxford Science Publications. The Clarendon Press, Oxford University Press, Oxford, 2004.


\bibitem{14}Brown, L.G.; Douglas, R. G.; Fillmore, P. A. Extensions of $C\sp*$-algebras and $K$-homology. Ann. of Math. (2) \textbf{105} (1977), no. 2, 265--324.


\bibitem{1}Busby, R.C., Double centralizers and extensions of $C\sp{\ast} $ -algebras. Trans. Amer. Math. Soc. \textbf{132} (1968) 79--99.

\bibitem{10}Conway, J.B., A course in functional analysis. Graduate Texts in Mathematics, 96. Springer-Verlag, New York, 1985.

\bibitem{11}Dieudonne, J., Foundations of modern analysis. Enlarged and corrected printing. Pure and Applied Mathematics, Vol. 10-I. Academic Press, New York-London, 1969.



\bibitem{17}Dugundji, J.,  An extension of Tietze's theorem. Pacific J. Math. \textbf{1}, (1951). 353--367.




\bibitem{15}Kasparov, G.G., The operator $K$ -functor and extensions of $C\sp{\ast} $ -algebras. (Russian) Izv. Akad. Nauk SSSR Ser. Mat. \textbf{44} (1980), no. 3, 571--636, 719.

\bibitem{7}MacLane, S., Categories for the working mathematician. Graduate Texts in Mathematics, Vol. 5. Springer-Verlag, New York-Berlin, 1971.

\bibitem{9}Murphy, G.J., $C\sp *$ -algebras and operator theory. Academic Press, Inc., Boston, MA, 1990.


\bibitem{6}Reid, G. A., Epimorphisms and surjectivity. Invent. Math. \textbf{9} (1969/1970), 295--307.
\bibitem{8}Rotman, J.J., An introduction to homological algebra. Second edition. Universitext. Springer, New York, 2009.




\bibitem{4}Rudin, W., Principles of mathematical analysis. Third edition. International Series in Pure and Applied Mathematics. McGraw-Hill Book Co., New York-Auckland-DCÂ¼sseldorf, 1976.






\bibitem{2}Yang, K.W., Completion of normed linear spaces. Proc. Amer. Math. Soc. \textbf{19} (1968), 801--806.

\end{thebibliography}
\end{document}